\documentclass{article}
\usepackage{PRIMEarxiv}

\RequirePackage{fix-cm}
\usepackage{graphicx}
\usepackage{amsmath,amssymb,amsfonts}
\usepackage{mathrsfs}
\usepackage{wasysym}
\usepackage{algorithm}
\usepackage{algorithmic}

\usepackage{enumitem}
\usepackage{bm,bbm}
\usepackage{multirow,booktabs,comment}
\usepackage{epstopdf}
\usepackage{xcolor}
\usepackage{soul}
\usepackage{subcaption}
\usepackage[colorlinks=true,urlcolor=blue,citecolor=red,anchorcolor=black,linkcolor=blue]{hyperref}

\newcommand{\bx}{\boldsymbol{x}}
\newcommand{\by}{\boldsymbol{y}}
\newcommand{\bX}{\boldsymbol{X}}

\newcommand{\bu}{\boldsymbol{u}}

\newcommand{\bXheta}{\boldsymbol{\theta}}

\begin{document}

\title{Predictive Moving Sample Method for Physics-Informed Neural Solvers of Time-Dependent PDEs}

\author{
  Beining Xu$^{\star}$, Bocheng Zhang$^{\star}$ \\
  The School of Information Science and Technology and the Institute of Mathematical Sciences\\
  ShanghaiTech University \\
  Shanghai 201210, China \\
  \texttt{\{xubn2024, zhangbch2023\}@shanghaitech.edu.cn} \\
   \And
  Haijun Yu \\
  SKLMS \& LSEC, Institute of Computational Mathematics and Scientific/Engineering Computing \\
  Academy of Mathematics and Systems Science, Chinese Academy of Sciences \\
  Beijing 100190, China \\
  School of Mathematical Sciences \\
  University of Chinese Academy of Sciences \\
  Beijing 100049, China \\
    \texttt{hyu@lsec.cc.ac.cn} \\
   \AND
   Zhao Zhang \\
    Research Center for Mathematics and Interdisciplinary Sciences\\
    Shandong University \\
    Qingdao, Shandong Province, 266237, China \\
    Frontiers Science Center for Nonlinear Expectations, Minister of Education \\
    Shandong University \\
    Qingdao, Shandong Province, 266237, China \\
    \texttt{zhaozhang@sdu.edu.cn} \\
   \And
   Jiayu Zhai$^{\dagger}$ \\
  The Institute of Mathematical Sciences \\
  ShanghaiTech University \\
  Shanghai 201210, China \\
  \texttt{zhaijy@shanghaitech.edu.cn}
}
\renewcommand{\thefootnote}{}
\footnotetext[1]{$^{\star}$ Co-first authors listed alphabetically.}
\footnotetext[2]{$^{\dagger}$ Corresponding author.}

\maketitle

\begin{abstract}
Time-dependent partial differential equations (PDEs) often develop sharp fronts, localized peaks, and other moving structures that occupy only a small portion of the space--time domain but dominate the approximation error. This makes fixed or uniformly sampled collocation strategies inefficient for physics-informed neural networks (PINNs), especially in high dimensions and over long-time prediction intervals. We propose the predictive moving sample method (PMSM), which builds on the moving sample method (MSM) in \cite{xu2026moving} by replacing its full time domain iterative training with a progressive time-stepping strategy and simplifying the velocity-field loss to further reduce the per-step cost. To improve practicality for long-time prediction, we further introduce the windowed-reset predictive moving sample method (WR-PMSM), which restricts extension training to an active time window and periodically resets the reference state, thereby reducing the growth of optimization cost while preserving global consistency through a final refinement stage. Across four representative benchmarks, PMSM consistently outperforms both standard PINNs and the original MSM under matched collocation budgets. These results suggest that transporting samples according to residual dynamics provides an effective and practical route to neural network solvers for time-dependent PDEs.
\keywords{Partial differential equations \and Physics-informed neural networks \and Transport equation \and Adaptive sampling \and Optimal transport}
\end{abstract}

\section{Introduction}

Partial differential equations (PDEs) provide a standard mathematical framework for describing various physical phenomena, such as fluid dynamics, heat conduction, and wave propagation. In recent years, physics-informed neural networks (PINNs) have developed rapidly. By directly embedding physical PDEs and initial or boundary conditions into the training process of neural networks and utilizing automatic differentiation in machine learning frameworks, PINNs conveniently train neural network solution approximations at the sampled collocation points in the computational domain \cite{raissi2019physics}. As a mesh-free method, PINNs avoid mesh construction and numerical differentiation and provide a flexible framework for solving both forward and inverse high-dimensional problems. However, this brings a new problem: when the equations involve both high dimensionality and moving singularities, the sampling strategy of the training collocation points becomes a key factor determining the training efficiency and solution quality.

Vanilla PINNs usually just use uniformly sampled collocation points most of which are not able to capture singularity characteristics of the solution, so that the neural network cannot learn these characteristics. Starting from \cite{wu2023rar}, adaptive sampling strategies were introduced to overcome this difficulty for neural network solvers. As a counterpart of adaptive mesh refinement in traditional numerical methods such as finite element method (FEM), most adaptive sampling methods are designed with ideas from adaptive mesh refinement. In \cite{tang2023das}, an adaptive sampling method were explained as a variance reduction technique through importance sampling, which can be considered as a starting point of theoretical study of adaptive sampling method. In \cite{tang2024adversarial}, adaptive sampling is considered as an optimal collocation distribution selection for integrating the differential equation, and the convergence of the method and the adaptivity of the samples are proved with optimal transport theory.

The aforementioned methods mainly deal with stationary equations, which means the singularities stay in local regions. However, a major difficulty in solving time-dependent PDEs is that these equations are not complex over the entire spatio-temporal domain and their local singularities move over time. Then the training collocation points need to capture the dynamics of the moving singularities, so that the neural network can focus on learning them efficiently, while spending less effort on relatively smoother regions. Regarding this characteristic of time-dependent PDEs, some works have studied the evolution of samples over time. For example, active learning and particle dynamics in the neural Galerkin direction \cite{bruna2024neural,Yuxiao_Coupling_2023}, and the normalized residual distribution at discrete time steps is treated as the invariant distribution of stochastic differential equations (SDEs), then these SDEs are used to generate adaptive samples. Similarly, \cite{WangHu2024NonuniformRandom,LiuZDYu2026} also use SDEs to sample moving collocation points. However, SDE samplers need to subdivide time steps into sub-intervals, which requires the numerical SDE scheme to generate transient samples converging to the residual distribution fast (see \cite{dobson2021using}). In \cite{gao2023active}, MCMC is used to regenerate adaptive samples at each time step, but this causes a huge computational cost and relies on the design of the MCMC to avoid high rejection rates, particularly for high singularities and high dimensionalities of the problem. In \cite{yang2024moving}, a moving sampling PINNs inspired by moving mesh PDEs is introduced, but it is significantly different from our approach. Our approach is based on moving sample method (MSM) \cite{xu2026moving}, which can be considered as the probabilistic version of the moving mesh method. To be specific, the dynamics of the adaptive samples are driven by a gradient ordinary differential equation (ODE) system whose potential is learned from a governing equation. In \cite{xu2026moving}, a probabilistic version of equidistribution property is both theoretically proved and numerically verified for MSM. 

Although MSM can efficiently capture the dynamics of singularities with the coupled gradient system, it still needs to learn the time-dependent potential over the entire time interval. This makes the learning for the ODE sampler expensive if the time interval is long and difficult for MSM to make long-time prediction. To address this we propose the predictive moving sample method (PMSM). Instead of training the potential on all $N_t$ time steps, PMSM introduces a progressive time-stepping strategy. To be specific, PMSM trains on only $N_t^{\text{init}}$ time steps initially, where $N_t^{\text{init}}$ is a relatively small integer. Then it extends on the time dimension one step at a time: using the current velocity field to predict temporary samples for the new time step, training the PINNs, updating the velocity field with the new residuals, and finally evolving the adaptive samples with the updated velocity field. This progressive strategy allows the method to iterate over longer time domains and improves the capacity to capture long-term solution trajectories. At the same time, when the time step is chosen small enough, we can assume the residual is approximately constant. Then the governing equation for the adaptive sample distribution is significantly simplified without the global integral-normalization term at each time step in the original MSM loss. The local velocity learning and simplified loss together reduce the overall computational cost and make long-time prediction possible. After the extension stage, the velocity networks will be dropped and only the adaptive samples over the entire time domain are collected. Then a final refinement stage trains the solution network on the full set of collected samples.

To make this method more practical for long-time prediction, we further develop the windowed-reset predictive moving sample method (WR-PMSM). Because the simplified velocity-field loss no longer depends on a global residual integral, we can restrict extension-stage PINNs training to an active time window and periodically reset the temporal reference point for the initial condition. At each reset, the current PINNs is frozen as a reference model that supplies the initial-condition supervision for the next window, while the velocity network continues to track singular regions on local time slices. After all samples have been collected, a global refinement stage restores consistency over the entire time domain. This strategy inhibits the growth of per-step training complexity with respect to the time-domain size, improving the usability of PMSM on problems with long evolution times.

The remaining part of this paper is organized as follows. Section~\ref{Preliminaries} introduces preliminaries, including moving mesh methods, PINNs, collocation distributions and residual-based adaptivity, probability transport and sample flows, HMC-based initial sampling, and a review of the MSM formulation and algorithm. Section~\ref{PMSM} presents the PMSM method design, covering the simplified velocity-field loss, the progressive time-extension strategy, and the WR-PMSM variant. Section~\ref{Numerical Experiments} reports numerical results. Section~\ref{conclusion} concludes the paper.

\section{Preliminaries}\label{Preliminaries}

This section sets the notation and mathematical background for PMSM. We first review the moving-mesh idea that motivates sample dynamics for MSM, then recall the residual formulation of PINNs for time-dependent PDEs, discuss how the choice of collocation distribution affects the empirical residual loss, summarize the probability-transport viewpoint used to connect moving samples with evolving residual distributions, describe the Hamiltonian Monte Carlo (HMC) initialization used to seed the adaptive trajectories, and finally provide a systematic review of the moving sample method (MSM) formulation and algorithm as the foundation on which PMSM builds.

\subsection{Moving mesh PDE methods (MMPDE)}

Moving mesh methods provide a classical paradigm for dynamic resource allocation in numerical PDEs. It lies the equidistribution principle \cite{boor1973good,white1979selection}: mesh density should be controlled by a positive monitor function $M(\bx)$ that measures local computational demand. In one dimension, the coordinate mapping $x = x(\xi)$ from a computational coordinate $\xi\in[0,1]$ to the physical domain satisfies
\begin{equation}
M(x(\xi))\,\frac{\partial x}{\partial \xi} = C,
\end{equation}
where $C$ is a constant independent of $\xi$, determined by boundary conditions so that the total number of mesh points remains fixed. A larger monitor value implies finer physical mesh spacing. Typical choices of $M$ rely on local error indicators, gradients, curvature, or residual information.

For time-dependent problems, MMPDE relocate grid points by continuously evolving the coordinate mapping in time. A common construction starts from a mesh-energy functional \cite{huang2010adaptive}:
\begin{equation}
    I[\boldsymbol{\xi}] = \frac{1}{2} \int_\Omega \sum_{i=1}^d \nabla \boldsymbol{\xi}_i^T M_i^{-1} \nabla \boldsymbol{\xi}_i d\bx,
\end{equation}
where each $M_i$ is a symmetric positive-definite monitor matrix. The gradient-flow relaxation
\begin{equation}
    \frac{\partial \boldsymbol{\xi}}{\partial t} = - \frac{P}{\tau} \frac{\delta I}{\delta \boldsymbol{\xi}},
\end{equation}
uses a parameter $\tau>0$ to control the response time of mesh motion and an operator $P$ to select a specific MMPDE form \cite{huang1994MMPDES}. Near the steady state of the energy, this yields
\begin{equation}
\frac{\partial \boldsymbol{\xi}_i}{\partial t}
= \frac{P}{\tau}\nabla\cdot(M_i^{-1}\nabla\boldsymbol{\xi}_i).
\end{equation}
In one dimension, a standard instance of MMPDE is given by
\begin{equation}
\frac{\partial x}{\partial t} = \frac{1}{\tau} \frac{\partial}{\partial \xi}\left(M \frac{\partial x}{\partial \xi}\right).
\end{equation}

The moving-mesh viewpoint separates two concerns: the monitor function identifies where computational resolution is needed, while the mesh dynamics describe how degrees of freedom should move. MSM \cite{xu2026moving} bears a similar spirit: the PDE residual plays the role of a monitor, while the adaptive collocation points move according to a learned velocity field rather than a prescribed mesh equation.

\subsection{Introduction to physics-informed neural networks}\label{subsec:pinns}

Physics-informed neural networks (PINNs) \cite{raissi2019physics} have emerged as a novel deep learning framework for solving partial differential equations. By embedding the governing physical laws into the loss function, PINNs enforce consistency between the network output and the underlying PDE constraints. This integration enables the method to leverage both data and physics, offering a flexible alternative to traditional numerical solvers. The general form of the problem is a time-dependent PDE in $\Omega \times [0,T]$:
\begin{equation}\label{eq:general_time_pde}
\left\{
\begin{aligned}
\mathcal{L}u(\bx,t) &= f(\bx,t), &&(\bx,t)\in\Omega\times[0,T],\\
\mathcal{I}u(\bx,0) &= u_0(\bx), &&\bx\in\Omega,\\
\mathcal{B}u(\bx,t) &= u_b(\bx,t), &&(\bx,t)\in\partial\Omega\times[0,T].
\end{aligned}
\right.
\end{equation}
Here $\mathcal{L}$ denotes the differential operator of the governing equation, while $\mathcal{I}$ and $\mathcal{B}$ denote the initial and boundary operators. 

In PINNs, the unknown solution $u(\mathbf{x},t)$ is approximated by a parameterized deep neural network, which is usually implemented as a fully connected feedforward neural network $u_{\boldsymbol{\theta}}(\mathbf{x},t)$ with $L$ layers. 
Let $\boldsymbol{\theta}$ denote the collection of all trainable parameters. Denoting the 
input as $\mathbf{h}^{(0)} = (\mathbf{x}, t)$, the network computes
\begin{gather*}
\mathbf{h}^{(l)} = \sigma\!\big(W^{(l)}\mathbf{h}^{(l-1)} + b^{(l)}\big),\quad l=1,\dots,L-1 ,\\
u_{\boldsymbol{\theta}}(\mathbf{x},t) = W^{(L)}\mathbf{h}^{(L-1)} + b^{(L)},
\end{gather*}
where $W^{(l)}$ and $b^{(l)}$ are the weight matrix and bias vector of layer $l$, and 
$\sigma(\cdot)$ is a nonlinear activation function. The interior residual is
\begin{equation}\label{eq:pde_residual_prelim}
r_{\bXheta}(\bx,t)=\mathcal{L}u_{\bXheta}(\bx,t)-f(\bx,t).
\end{equation}
For notational convenience, we write $r_t(\bx)$ as a shorthand for $r_{\bXheta}(\bx,t)$ at a fixed time $t$ throughout the paper. The initial and boundary residuals are
\begin{equation}\label{eq:ic_bc_residual_prelim}
r_{\bXheta}^{0}(\bx)=\mathcal{I}u_{\bXheta}(\bx,0)-u_0(\bx),
\qquad
r_{\bXheta}^{b}(\bx,t)=\mathcal{B}u_{\bXheta}(\bx,t)-u_b(\bx,t).
\end{equation}

The overall PINNs objective is a weighted residual minimization problem. With positive weights $\beta_0$ and $\beta_b$, a standard continuous loss takes the form
\begin{align}
\mathcal{J}_{\mathrm{PINNs}}(\bXheta)
=&
\int_0^T\int_{\Omega}|r_{\bXheta}(\bx,t)|^2\,d\bx dt
{}+\beta_0\int_{\Omega}|r_{\bXheta}^{0}(\bx)|^2\,d\bx \notag\\
&+\beta_b\int_0^T\int_{\partial\Omega}|r_{\bXheta}^{b}(\bx,t)|^2\,dS\,dt .
\label{eq:pinn_loss_prelim}
\end{align}
In practice, the integrals in \eqref{eq:pinn_loss_prelim} are replaced by Monte Carlo averages over interior, initial, and boundary point sets. Thus the empirical training objective depends not only on the network architecture and optimizer, but also on the probability measures from which the collocation points are drawn. This dependence is crucial in the present work, because our goal is to replace a fixed interior sampling measure by a time-dependent adaptive measure.

\subsection{Collocation distributions and residual-based adaptivity}

The distribution of residual points determines which regions dominate the empirical PDE loss. If the interior points are sampled uniformly from $\Omega\times[0,T]$, the empirical loss approximates an average residual over the whole domain. Such averaging is efficient for smooth solutions but can be ineffective when the solution contains localized fronts, peaks, or high-gradient layers. In those cases, the complex region may occupy a very small volume, so its contribution to the global average can be diluted even when its local error is large. This sampling effect is one of the reasons why standard PINNs can fail or train slowly on convection-dominated and strongly time-dependent problems \cite{krishnapriyan2021characterizing,wang2021gradient}.

Residual-based adaptive sampling addresses this issue by changing the point distribution according to information obtained during training. Existing strategies include residual-based refinement, adaptive residual distributions, active learning, failure-informed sampling, generative sampling, and optimal-transport-based sampling \cite{wu2023rar,gao2023active,gao2023failure,tang2023das,tang2024adversarial}. A common abstraction behind these methods is to construct a residual-induced density \cite{xu2026moving}, for example
\begin{equation}\label{eq:residual_induced_density}
\rho_t(\bx)=\frac{|r_t(\bx)|^2+\varepsilon}{\int_\Omega (|r_t(\by)|^2+\varepsilon)\,d\by},
\end{equation}
where $r_t(\bx)$ is the PDE residual at time $t$ (as defined in Section~\ref{subsec:pinns}) and $\varepsilon>0$ provides a uniform background sampling while preventing degeneracy. Sampling from $\rho_t$ allocates more points to high-residual regions and therefore increases their influence in the empirical loss.

For time-dependent PDEs, however, residual-based adaptivity should account for temporal motion as well as spatial concentration. A high-residual region at time $t$ may move, deform, or disappear at a later time. Repeatedly resampling independent point sets at each time slice can respond to the current residual, but it does not preserve the trajectories of previously informative points, nor does it exploit the continuity of the solution structures between adjacent time steps. MSM \cite{xu2026moving} addresses this by treating adaptive collocation as the evolution of a sample distribution: adaptive points are viewed as particles transported by a velocity field, continuously following the residual structures rather than being regenerated at each time slice. We adopt the same perspective in this work.

\subsection{Probability transport and sample flows}\label{subsec:prob_transport}

In MSM, moving samples are described as a particle system transported by a velocity field \cite{xu2026moving}. Let $\Phi_t:\Omega\to\Omega$ be a flow map and let $\mathbf{V}(\bx, t)$ be a velocity field satisfying
\begin{equation}\label{eq:flow_ode_prelim}
\frac{\partial \Phi_t}{\partial t}(\bx)
=
\mathbf{V}(\Phi_t(\bx), t),
\qquad
\Phi_0(\bx)=\bx .
\end{equation}
If the initial sample density is $\rho_0$ and $\rho_t$ denotes the density induced by the pushforward of $\rho_0$ through $\Phi_t$, then $\rho_t$ satisfies the continuity equation
\begin{equation}\label{eq:continuity_prelim}
\frac{\partial \rho_t}{\partial t}
+\nabla\cdot(\rho_t\mathbf{V}_t)=0 .
\end{equation}
Conversely, under suitable regularity assumptions, solutions of \eqref{eq:continuity_prelim} describe the density evolution generated by the particle flow \eqref{eq:flow_ode_prelim} \cite{ambrosio2008gradient}. This equation is the probability analogue of conservation of mass: the probability density changes because samples are transported by the velocity field. It serves as the theoretical foundation for the sample evolution in both MSM and the proposed PMSM.

\subsection{Initial adaptive sampling by Hamiltonian Monte Carlo}\label{subsec:hmc_prelim}

Before the sample flow is evolved, the initial adaptive particles should already resolve the important structures of the initial condition. Uniform initialization is simple, but it can place too few particles near narrow peaks, transition layers, or high-gradient fronts at $t=0$. We therefore generate the initial point set by sampling from an initial-condition-induced density. For a generic nonnegative initial profile, we use
\begin{equation}\label{eq:initial_mcmc_density}
\pi_0(\bx)
=
\frac{|u_0(\bx)|+\varepsilon}
{\int_{\Omega}(|u_0(\by)|+\varepsilon)\,d\by},
\end{equation}
where $\varepsilon>0$ is a small stabilizing constant. For Burgers-type traveling-front problems, sampling directly from $u_0$ would overemphasize the plateau of the sigmoid initial condition. In these cases, we instead use the gradient-energy monitor
\begin{equation}\label{eq:initial_burgers_mcmc_density}
\pi_0^{\mathrm{B}}(\bx) = \frac{\|\nabla u_0(\bx)\|_2^2+\varepsilon}{\int_{\Omega}(\|\nabla u_0(\by)\|_2^2+\varepsilon)\,d\by},
\end{equation}
which concentrates samples around the initial wave front. Both densities allocate more points to regions where the solution has stronger localized structure.

We approximate samples from \eqref{eq:initial_mcmc_density} or \eqref{eq:initial_burgers_mcmc_density} by Hamiltonian Monte Carlo (HMC), an MCMC method \cite{duane1987hybrid}. Its proposal is corrected by the Metropolis--Hastings accept-reject rule \cite{metropolis1953equation,hastings1970monte}. The target log probability is
\begin{equation}\label{eq:initial_log_prob}
\ell_0(\bx)=
\begin{cases}
\log(|u_0(\bx)|+\varepsilon), & \bx\in\Omega\quad\text{for generic initial profiles},\\
\log(\|\nabla u_0(\bx)\|_2^2+\varepsilon), & \bx\in\Omega\quad\text{for Burgers-type equations},\\
-\infty, & \bx\notin\Omega .
\end{cases}
\end{equation}
Given a current state $\bx$ and an auxiliary momentum $\boldsymbol{p}$, HMC simulates leapfrog steps for the Hamiltonian
\begin{equation}\label{eq:hmc_hamiltonian}
H(\bx,\boldsymbol{p})=-\log \pi_0(\bx)+\frac{1}{2}\|\boldsymbol{p}\|_2^2,
\end{equation}
and accepts the proposal $(\bx',\boldsymbol{p}')$ with probability
\begin{equation}\label{eq:hmc_acceptance}
\alpha(\bx,\bx')
=
\min\{1,\exp[-H(\bx',\boldsymbol{p}')+H(\bx,\boldsymbol{p})]\}.
\end{equation}
In the implementation, multiple chains are initialized in parallel, a burn-in phase is used to adapt the step size toward a target acceptance probability, and the resulting samples are shuffled and filtered to keep only points inside the computational domain. The accepted MCMC samples provide two point sets: the first part is used in $\mathcal{C}_0$ for the initial loss, while the remaining part is used as the seed set $\mathbf{Z}_0$ for the moving adaptive trajectories. We additionally mix a small number of uniformly sampled initial points into $\mathcal{C}_0$ to keep the initial-condition penalty globally anchored.

\subsection{Moving sample method (MSM)}\label{subsec:MSM}

The core idea of moving sample method (MSM) \cite{xu2026moving} is to replace the fixed uniform sampling in standard PINNs with a time-evolving collocation measure $\nu_t$, where $\nu_t$ is determined by the velocity field network trained by the transport equation \eqref{eq:continuity_prelim}. By setting an appropriate target distribution, sampling points can be clustered in regions with high singularity.

Based on transport equation, MSM typically chooses $r_t^2$ as the unnormalized target sampling distribution. This choice is consistent with the squared-residual objective of PINNs: regions that contribute more to the PDE loss should receive higher collocation density.

MSM employs two fully connected neural networks: a solution network $u_{\bXheta}(\bx,t)$ and a scalar potential network $\psi_{\boldsymbol{\omega}}(\bx,t)$, with the velocity field obtained as $\mathbf{V}_{\boldsymbol{\omega}}=\nabla\psi_{\boldsymbol{\omega}}$. The solution network is trained with the standard PINNs residual loss, except that the measure of collocation points is replaced by the dynamic measure $\nu_t$, $\nu_0$ and $\nu_b$. The loss of $u_{\bXheta}$ is given by
\begin{align}
    \mathcal{J}_u(\bXheta) = & \int_0^T\int_\Omega |\mathcal{L}u_{\bXheta} - f|^2 \,\nu_t(d\bx)dt + \lambda_0 \int_{\Omega} |\mathbf{\mathcal{I}}u_{\bXheta} - \bu_0|^2 \nu_0(d\bx) \notag\\
    & + \lambda_b \int_0^T\int_{\partial \Omega} |\mathcal{B}u_{\bXheta} - u_b|^2 \,\nu_b(d\bx)dt, \label{loss_u}
\end{align}
where $\lambda_0, \lambda_b > 0$ are balancing weights. 

For the velocity network $\mathbf{V}_{\boldsymbol{\omega}}$, applying the logarithmic form of the continuity equation to the target density yields the constraint that the velocity field $\mathbf{V}_t$ should satisfy:
\begin{equation}\label{eq:velocity_constraint}
2\frac{\partial r_t}{\partial t} + 2\nabla r_t\cdot\mathbf{V}_t + r_t\,\nabla\cdot\mathbf{V}_t - \frac{r_t}{\int_{\Omega} r_t^2\,d\bx}\,G_t = 0, \qquad G_t:=\frac{d}{dt}\int_{\Omega} r_t^2\,d\bx .
\end{equation}

Equation~\eqref{eq:velocity_constraint} alone does not determine $\mathbf{V}_t$ uniquely, as it constrains only the divergence of the velocity field. By the Helmholtz decomposition \cite{arfken2011mathematical}, the solenoidal component does not affect density evolution and thus constitutes a degree of freedom. MSM therefore restricts $\mathbf{V}_t$ to a pure gradient form, parameterized through a scalar potential network $\psi_{\boldsymbol{\omega}}$, i.e. $\mathbf{V}_{\boldsymbol{\omega}}=\nabla\psi_{\boldsymbol{\omega}}$, which reduces the output dimensionality of the network $\psi_{\boldsymbol{\omega}}$.

Substituting $\mathbf{V}_{\boldsymbol{\omega}}=\nabla\psi_{\boldsymbol{\omega}}$ into~\eqref{eq:velocity_constraint} and converting to a least-squares objective gives the loss of velocity network $\psi_{\boldsymbol{\omega}}$:
\begin{align}
\mathcal{J}_{\mathbf{V}}^{\mathrm{MSM}}({\boldsymbol{\omega}}) = \int_0^T\int_\Omega \Big( 2\frac{\partial r_t}{\partial t} + 2\nabla r_t\cdot\nabla\psi_{\boldsymbol{\omega}} + r_t\Delta\psi_{\boldsymbol{\omega}} - \frac{r_t}{\int_\Omega r_t^2\,d\bx}G_t \Big)^2 \,d\bx\,dt , \label{loss_v_MSM}
\end{align}
where $G_t$ is approximated by $$
G_t = \frac{d}{dt}\int_{\Omega} r_t^2\,d\bx \approx \frac{\int_{\Omega} r_{t_{i+1}}^2 - r_{t_i}^2\,d\bx}{t_{i+1}-t_i}.
$$

The MSM algorithm can be summarized as three stages: pretrain the network $u_{\bXheta}$ on uniform sampling points to get a velocity field $\mathbf{V}_t$ firstly; then, for several iterations, add new collocation points by using the current velocity field to evolve the samples over all $N_t$ time steps and retrain both networks; finally, after all iterations, perform a final refinement training of the solution network on the collected samples.

MSM provides an elegant probability-transport framework for adaptive sampling in PINNs for solving time-dependent PDEs. However, it can be inefficient when dealing with some dynamical systems with long evolution times, for the following two reasons:

(1) The loss function of the velocity field $\mathbf{V}_t$ requires significant computational resources. Since $G_t$ contains the time derivative with respect to the integral, time slices need to be passed in pairs during computation. This affects computational efficiency to some extent. Furthermore, the original method is also difficult to solve the evolution of dynamical systems over long time periods. When time $T$ is large (corresponding to a large number of time steps $N_t$), the collocation set $\mathcal{C} = \{(\bx_{\Omega}^{(i)}, t^{(i)})\}_{i=1}^{N\times N_t}$ grows linearly with $N_t$, causing memory consumption and computation time to increase steadily per iteration.

(2) Re-evolving the adaptive samples on all time steps from the initial points at every iteration reduces efficiency for long-time prediction. MSM resamples $\mathbf{Z}_0$ from $\pi_0$ at each iteration, then evolves samples over all $N_t$ time steps with the newly trained velocity field. This is reasonable for problems over short time domains---it allows each round to re-evolve the adaptive trajectories using an improved solution estimate. In long-time scenarios, however, information from the initial time $t=0$ has limited influence on residuals near $t=T$, so the benefit of repeatedly evolving all time steps from scratch diminishes while the computational cost continues to grow. Moreover, $\mathcal{J}_{\mathbf{V}}^{\mathrm{MSM}}$ involves a double integral over the full spatial and temporal domains; recomputing this global quantity at every iteration adds considerable overhead.

\section{Predictive moving sample method (PMSM)}\label{PMSM}

Building upon the probability transport framework of Moving sample method (MSM), we introduce an improved method named Predictive moving sample method (PMSM). The main modifications can be summarized in the following two aspects. First, the extension stage adopts a progressive time-stepping strategy, which focuses on the local evolution of adaptive collocation points without requiring a fully trained solution and velocity network over the entire time domain at every iteration. Second, and in connection with this, the velocity field loss is simplified: the time integral is restricted to the current time slice, and the global normalization term is dropped under a small timestep approximation.

\subsection{Loss simplification for velocity field}\label{subsec:pmsm_loss}

For $\mathcal{J}_{\mathbf{V}}^{\mathrm{MSM}}$ in equation \eqref{loss_v_MSM}, the $G_t$ involves computationally demanding operations, which significantly affect the efficiency of the algorithm. Since $G_t$ involves taking the time derivative of a spatial integral of $r_t^2$, its evaluation requires a discrete approximation, which makes the loss of $\mathbf{V}_t$ more complex. Computing the loss at a given time $t_1$ requires access to the residual at the next time step $t_2$, so the loss input $(x, t)$ is replaced by $(x, (t_1, t_2))$ in MSM and evaluates the residuals at both times. Although we currently adopt a fixed time step $\Delta t$, this dependency still introduces a computational burden, which becomes particularly pronounced in long-time simulations.

Motivated by the above issue, we seek to simplify the term $G_t$. A natural observation is that the residuals at two adjacent time typically do not vary significantly, an effect that becomes even more pronounced when the time step $\Delta t$ is small enough. This property can be further strengthened through appropriate algorithmic design, enabling us to neglect this term under suitable conditions. Therefore the computation of the full time-domain integral in $\mathcal{J}_{\mathbf{V}}^{\mathrm{MSM}}$~\eqref{loss_v_MSM} reduces to computation on the current time slice.

Concretely, PMSM assumes that the global normalization term $\frac{r_t}{\int_\Omega r_t^2 d\bx} G_t$ (where $G_t = \frac{d}{dt} \int_\Omega r_t^2 d\bx$) is small enough, since under the small timestep hypothesis the variation of the normalization constant is small compared with the local residual variation. Compared with the original MSM, PMSM no longer requires numerical evaluation of the derivative terms with respect to the integral in the $\mathcal{J}_{\mathbf{V}}$, which significantly reduces the training time.

In summary, the simplified PMSM velocity-field loss is
\begin{align}
    \mathcal{J}_{\mathbf{V}}^{\text{PMSM}}({\boldsymbol{\omega}}) &= \int_\Omega \Big(2\frac{\partial}{\partial t}r_t + 2\nabla r_t \cdot \mathbf{V}_{\boldsymbol{\omega}} + r_t \nabla\cdot \mathbf{V}_{\boldsymbol{\omega}}\Big)^2 \,d\bx \notag\\
    &= \int_\Omega \Big(2\frac{\partial}{\partial t}r_t + 2\nabla r_t \cdot \nabla \psi_{\boldsymbol{\omega}} + r_t \Delta \psi_{\boldsymbol{\omega}}\Big)^2 \,d\bx . \label{loss_v_PMSM}
\end{align}

\subsection{Progressive time extension strategy}\label{subsec:pmsm_extension}

The core innovation of PMSM is a progressive time-extension strategy that replaces MSM's full-time-domain synchronous iteration with step-by-step time marching. For boundary and initial conditions, the collocation points are fixed spatial point sets repeated over all time steps, following the same approach as MSM. For the adaptive collocation points used in the PDE loss, we propose a new training procedure, detailed below.

The first stage is the initial stage for $u_{\bXheta}$ and $\mathbf{V}_{\boldsymbol{\omega}}$. PMSM trains on only the first $N_t^{\text{init}}$ time steps rather than total time steps. The initial time grid is given by
\begin{equation}\label{eq:initial_time_grid}
\mathcal{T}_{\mathrm{global}}^{0}=\{t_0,t_1,\dots,t_{N_t^{\mathrm{init}}-1}\}.
\end{equation}
On this initial block, the solution network $u_{\bXheta}$ is trained first, followed by the potential network $\psi_{\omega}$. Finally, the adaptive seeds $\mathbf{Z}_0$ are pushed forward by the learned velocity field to form the initial adaptive trajectories over $\mathcal{T}_{\mathrm{global}}^{0}$. This stage essentially corresponds to one iteration of the MSM algorithm, except that it is carried out fewer shorter time steps. This stage provides PMSM with its initial solution model, initial adaptive sample set on $\mathcal{T}_{\mathrm{global}}^{0}$, and initial velocity field.

The second stage is the extension stage, which iterates several times, each time extending the computational time domain of the network by a certain length. In the experiments presented in Section~\ref{Numerical Experiments}, we prolong the time domain by exactly one time step per iteration. For each newly introduced time step $t_i$, PMSM follows a predict--train--update cycle:
 
\begin{enumerate}
    \item Temporary sampling. The solution network $u_{\bXheta}$ and velocity field $\mathbf{V}_{\boldsymbol{\omega}}$ are already available on $[t_0, t_{i-1}]$. Using the velocity field from the previous iteration, we can advance the existing adaptive trajectories from $t_{i-1}$ to $t_i$, producing temporary sampling points at the newly added time slice $t_i$.
    \item Training for solution network $u_{\bXheta}$. These temporary samples are combined with the existing adaptive collocation set as the training set for this round. The solution network $u_{\bXheta}$ is trained on $[t_0, t_i]$, maintaining global validity of the solution while obtaining residual information at $t_i$.
    \item VNN update and samples finalization. The velocity network $\psi_{\boldsymbol{\omega}}$ is trained using the residuals obtained on $[t_{i-1}, t_i]$. The updated velocity field is then used to formally advance the adaptive trajectories to $t_i$, yielding the final adaptive sampling points for that time step, which are added to the adaptive collocation set.
\end{enumerate}

In the implementation, PMSM maintains two time axes: $\mathcal{T}_{\mathrm{global}}$ records all predicted time steps and is used for trajectory rollout, evaluation, and final global refinement; $\mathcal{T}_{\mathrm{train}} \subseteq \mathcal{T}_{\mathrm{global}}$ is the active training window, on which the PINNs collocation points are built during the extension stage. In the basic PMSM, $\mathcal{T}_{\mathrm{train}} = \mathcal{T}_{\mathrm{global}}$, i.e., all extended time steps are used. The WR-PMSM variant in Section~\ref{subsec:wr_pmsm} will further reduce the extension-stage cost by limiting the size of $\mathcal{T}_{\mathrm{train}}$.

During the extension stage, we deliberately limit the extent of solution network training, using fewer epochs in practice. This is because the goal of solution-network updates at this stage is to produce a residual field that is sufficiently informative for training the velocity network, rather than to fully converge the solution at every newly added time step. After $u_{\bXheta}$ update, the velocity network is trained on the recent time slice, and the updated velocity field is used to formally advance the adaptive trajectories. The ultimate purpose of the extension stage is to collect adaptive samples across the entire time domain, ensuring that adaptive points correspond to high-residual regions, while keeping the intermediate optimization local and computationally affordable.

In this stage, we assume that when the time interval $\Delta t$ is sufficiently small, the variation of the velocity field between two consecutive time is limited. Under this assumption, the temporary sampling points generated are meaningful, as they can partially cover regions with strong singularities. This design also enables us to iterate the algorithm using loss \ref{loss_v_PMSM} for long time, since the temporary sampling points at $t_i$ reduce the discrepancy between the residual integrals at $t_{i-1}$ and $t_i$, there by diminishing the error introduced by the assumption underlying our simplified loss formulation.

The third stage after initial and extension is final global refinement. Once all extension rounds have finished, we reconstruct collocation points over the complete time grid $\mathcal{T}_{\mathrm{global}}$. The final PDE training set contains both the uniform spatial points repeated over all global time steps and the adaptive trajectories collected during the extension stage. Boundary samples are also rebuilt over the entire time domain. The solution network is then trained once more on this complete set, using a larger final training budget than an individual extension round.

Algorithm~\ref{alg:basePMSM} summarizes the base PMSM procedure at the step of method design.

\begin{algorithm}
\caption{Base PMSM training procedure\label{alg:basePMSM}}
\begin{algorithmic}[1]
    \REQUIRE Initial networks $u_{\bXheta}, \mathbf{V}_{\boldsymbol{\omega}}$, the initial uniform collocation points for PDE loss $S_u = \{(\bx^{(i)}_{\Omega}, t^{(i)})\}_{i=1}^{N_u \times N_t^{\text{init}}}$, initial training set $S_0 = \{\bx_0^{(i)}\}_{i=1}^{N_0}$, boundary training set $S_{bdry} = \{(\bx_{\partial \Omega}^{(i)}, t^{(i)})\}_{i=1}^{N_b\times N_t^{\text{init}}}$, the number of training iteration $K_{ext}$, training epochs $K_u, K_w$ for $u_{\bXheta}, \mathbf{V}_{\boldsymbol{\omega}}$ in each iteration.
    \STATE Train $u_{\bXheta}$ and $\mathbf{V}_{\boldsymbol{\omega}}$ on $S_u, S_0, S_{bdry}$ using PMSM.
    \STATE Draw initial adaptive seeds $\mathbf{Z}_0$ from $\pi_0$ by HMC with a Metropolis--Hastings correction, then evolve them to $\Phi_t$ for $t = 1, 2, \dots, N_t^{\text{init}} - 1$, and let $S_{adaptive}=\{(\bx_{adaptive}^{(i)}, t^{(i)})\}_{i=1}^{N\times N_t^{\text{init}}}$.
    \FOR{$j = 1, 2, \dots, K_{ext}$}
        \STATE Evolve temporary adaptive samples $\bX'_{N_t^{\text{init}} - 1 + j}$ based on $\mathbf{V}_{\boldsymbol{\omega}}$, concatenate $\bX'_{N_t^{\text{init}} - 1 + j}$ with the existing adaptive points to obtain (temporary) $S_{adaptive}=\{(\bx_{adaptive}^{(i)}, t^{(i)})\}_{i=1}^{N\times (N_t^{\text{init}}+j)}$. Add the same boundary sampling points at $t=(N_t^{\text{init}}-1+j)\Delta t$, update $S_{bdry} = \{(\bx_{\partial \Omega}^{(i)}, t^{(i)})\}_{i=1}^{N_b\times (N_t^{\text{init}}+j)}$.
        \STATE Update model $u_{\bXheta}$ for $K_u$ epochs by descending the stochastic gradient of $\mathcal{J}_u$ \eqref{loss_u} on $S_u\cup S_{adaptive}, S_0, S_{bdry}$.
        \STATE Update model $\mathbf{V}_{\boldsymbol{\omega}}$ for $K_w$ epochs by descending the stochastic gradient of $\mathcal{J}_{\mathbf{V}}$ \eqref{loss_v_PMSM} on $S_u$.
        \STATE Evolve adaptive samples $\bX_{N_t^{\text{init}} - 1 + j}$ based on new $\mathbf{V}_{\boldsymbol{\omega}}$, concatenate $\bX_{N_t^{\text{init}} - 1 + j}$ with the existing adaptive points to obtain adaptive set $S_{adaptive}=\{(\bx_{adaptive}^{(i)}, t^{(i)})\}_{i=1}^{N\times (N_t^{\text{init}}+j)}$.
    \ENDFOR
    \STATE Train $u_{\bXheta}$ for $K_{final}$ epochs on $S_u \cup S_{adaptive}, S_0, S_{bdry}$ (final refinement).
    \ENSURE $u_{\bXheta}, S_{adaptive}$.
\end{algorithmic}
\end{algorithm}

\subsection{Windowed-reset strategy for long time prediction}\label{subsec:wr_pmsm}

For long time prediction, the basic PMSM appends only one new time step per extension round, but each round trains on the entire $\mathcal{T}_{\mathrm{global}}$ ($\mathcal{T}_{\mathrm{train}} = \mathcal{T}_{\mathrm{global}}$, see Section~\ref{subsec:pmsm_extension}), so the computational cost still accumulates as $\mathcal{T}_{\mathrm{global}}$ grows. To eliminate this growth, we introduce the WR-PMSM variant, which restricts extension stage training to a fixed size active window.

In WR-PMSM, the two time axes $\mathcal{T}_{\mathrm{global}}$ and $\mathcal{T}_{\mathrm{train}}$ introduced in Section~\ref{subsec:pmsm_extension} are explicitly separated: collocation points during the extension stage are rebuilt only on $\mathcal{T}_{\mathrm{train}}$, while $\mathcal{T}_{\mathrm{global}}$ is used only for trajectory rollout and final evaluation. We adopt a fixed size window policy: only the most recent $W$ time steps are retained.

When a window reset is triggered, the current solution network $u_{\bXheta}$ is cloned as a reference model $\widehat{u}$. In the next extension stage, the initial loss switches from the analytic initial condition $u_0(\bx)$ to the reference model's output at the window start, $\widehat{u}(\bx, t_{\mathrm{start}})$, where $t_{\mathrm{start}} = \min \mathcal{T}_{\mathrm{train}}$. This means that as $\mathcal{T}_{\mathrm{train}}$ advances and discards earlier time slices, the solution network is no longer required to refer directly back to the analytic initial condition at $t=0$. Instead, it is required to remain consistent with the previously learned solution at the entrance of the new window. This design preserves temporal continuity while enabling local training within each window.

After all extension rounds, the collocation points are rebuilt on the full $\mathcal{T}_{\mathrm{global}}$, and the solution network is trained with the analytic initial condition to restore full time-domain consistency.

Algorithm~\ref{alg:wr_pmsm} gives WR-PMSM as used in the long-time experiments.

\begin{algorithm}
\caption{WR-PMSM for long-time prediction\label{alg:wr_pmsm}}
\begin{algorithmic}[1]
    \REQUIRE Solution network $u_{\bXheta}$, velocity network $\mathbf{V}_{\boldsymbol{\omega}}$,
             the initial uniform collocation points for PDE loss $S_u = \{(\bx^{(i)}_{\Omega}, t^{(i)})\}_{i=1}^{N_u \times N_t^{\text{init}}}$,
             initial training set $S_0 = \{\bx_0^{(i)}\}_{i=1}^{N_0}$,
             boundary training set $S_{bdry} = \{(\bx_{\partial \Omega}^{(i)}, t^{(i)})\}_{i=1}^{N_b\times N_t^{\text{init}}}$,
             initial adaptive seeds $\mathbf{Z}_0$ drawn from $\pi_0$,
             initial time grid $\mathcal{T}_{\mathrm{global}} = \{t_0, t_1, \dots, t_{N_t^{\text{init}}-1}\}$,
             window size $W$, extension rounds $K_{ext}$,
             training epochs $K_u, K_w$ for $u_{\bXheta}, \mathbf{V}_{\boldsymbol{\omega}}$ in each round.
    \STATE Set $\mathcal{T}_{\mathrm{train}} = \mathcal{T}_{\mathrm{global}}$.
    \STATE Train $u_{\bXheta}$ for $K_u$ epochs on $S_u, S_0, S_{bdry}$ by descending $\mathcal{J}_u$ \eqref{loss_u}.
    \STATE Train $\mathbf{V}_{\boldsymbol{\omega}}$ for $K_w$ epochs on $S_u$ by descending $\mathcal{J}_{\mathbf{V}}$ \eqref{loss_v_PMSM}.
    \STATE Evolve $\mathbf{Z}_0$ to $\Phi_t$ for $t = t_1, \dots, t_{N_t^{\text{init}}-1}$ and construct $S_{adaptive} = \{(\bx_{adaptive}^{(i)}, t^{(i)})\}_{i=1}^{N \times N_t^{\text{init}}}$ on $\mathcal{T}_{\mathrm{global}}$.
    \FOR{$j = 1, 2, \dots, K_{ext}$}
        \STATE Append new time step $t_{N_t^{\text{init}}-1+j}$ to $\mathcal{T}_{\mathrm{global}}$ and $\mathcal{T}_{\mathrm{train}}$.
        \STATE Drop the earliest time step if $|\mathcal{T}_{\mathrm{train}}| > W$.
        \STATE Build temporary samples $\bX'_{N_t^{\text{init}}-1+j}$ at the new step from the current $\mathbf{V}_{\boldsymbol{\omega}}$ and build $S_{adaptive}^\text{temporary} = S_{adaptive} \cup \bX'_{N_t^{\text{init}}-1+j}$.
        \IF{$\min\mathcal{T}_{\mathrm{train}}$ has changed}
            \STATE Clone $u_{\bXheta}$ as a frozen reference model $\widehat{u}$.
            \STATE In $\mathcal{J}_u$, replace the analytic initial condition $u_0(\bx)$ with $\widehat{u}(\bx, t_{\mathrm{start}})$ at $t_{\mathrm{start}} = \min\mathcal{T}_{\mathrm{train}}$, applied to the initial loss on $\mathcal{T}_{\mathrm{train}}$.
        \ENDIF
        \STATE Rebuild $S_u, S_{bdry}$ on $\mathcal{T}_{\mathrm{train}}$; update $S_{adaptive}^\text{temporary}$ to retain only the trajectories on $\mathcal{T}_{\mathrm{train}}$.
        \STATE Update $u_{\bXheta}$ for $K_u$ epochs on $S_u \cup S_{adaptive}^\text{temporary}, S_0, S_{bdry}$ by descending $\mathcal{J}_u$ \eqref{loss_u}.
        \STATE Update $\mathbf{V}_{\boldsymbol{\omega}}$ for $K_w$ epochs on $S_u$ by descending $\mathcal{J}_{\mathbf{V}}$ \eqref{loss_v_PMSM}.
        \STATE Evolve $\bX_{N_t^{\text{init}}-1+j}$ from the updated $\mathbf{V}_{\boldsymbol{\omega}}$ and concatenate to $S_{adaptive}$ on $\mathcal{T}_{\mathrm{global}}$.
    \ENDFOR
    \STATE Rebuild $S_u, S_{adaptive}, S_{bdry}$ on the full $\mathcal{T}_{\mathrm{global}}$.
    \STATE Train $u_{\bXheta}$ for $K_{final}$ epochs on $S_u \cup S_{adaptive}, S_0, S_{bdry}$ with analytic initial condition $u_0(\bx)$ (final refinement).
    \ENSURE $u_{\bXheta}$ and $S_{adaptive}$ on $\mathcal{T}_{\mathrm{global}}$.
\end{algorithmic}
\end{algorithm}

\section{Numerical Experiments}\label{Numerical Experiments}

In this section, we evaluate PMSM on four representative time-dependent PDEs. The first benchmark is a 2D Burgers' problem with a narrow traveling front. The second is a 2D parabolic problem whose Gaussian peak moves along a curved trajectory while its amplitude grows and saturates. The third is a 3D Fokker--Planck problem whose probability density remains concentrated near a moving shell. And the fourth is a 6D variable-speed Burgers' problem, which combines high-dimensional sparsity with a decelerating front. Together these benchmarks test whether moving samples can track localized structures across different geometries and levels of difficulty.

\subsection{Experimental setup}

In terms of network architecture, the solution network for PMSM, MSM and PINNs adopts the same fully connected feedforward structure with three hidden layers. Each hidden layer has 64 neurons. The layout from the input layer to the output layer is $[d+1, 64, 64, 64, 1]$, where $d$ is the spatial dimension, with an additional time dimension. For the velocity potential network used in both PMSM and MSM, a single hidden layer configuration of $[d+1, 256, 1]$ is employed. Experiments indicate that this shallow network already possesses sufficient expressive capacity to represent the velocity field, while also offering improved training stability and reduced optimization cost.

The main accuracy metrics are the relative $L^2$ error
\[
\mathcal{E}_{L_2} = \frac{\| u_{\mathrm{pred}} - u \|_2}{\| u \|_2},
\]
and the absolute $L^\infty$ error
\[
\mathcal{E}_{L_\infty} = \max_{\bx,t} |u_{\mathrm{pred}}(\bx,t) - u(\bx,t)|,
\]
where $u_{\mathrm{pred}}$ denotes the predicted solution from the network, and $u$ denotes the analytical solution.

For all four experiments, the initial time is taken as $t_0=0$, and the time step is $\Delta t = 0.05$. The initial rollout $N_t^{\text{init}} = 2$ for Experiments~\ref{subsec:2d_burgers}, \ref{subsec:6d_burgers} and \ref{subsec:2d_parabolic}, and $N_t^{\text{init}} = 6$ for \ref{subsec:3d-fokker-planck}.  Since each extension round appends one new time step, the number of extension iterations is computed as $$ K_{ext} = \frac{T - N_t^{\text{init}} \times \Delta t - t_0}{\Delta t} + 1 .$$

When comparing different methods, we follow the training strategy of PMSM by partitioning the temporal dimension into uniform intervals of size $\Delta t$, and adopting the same number of collocation points at each time step. Specifically, PMSM and MSM each employ $N$ adaptive points and $N_u$ uniform PDE points, while PINNs uses $N+N_u$ fixed points. MSM is also configured with the same network architectures, and total epochs as PMSM, and the only difference lies in the training procedure (iterative training over the full time domain for MSM versus progressive time extension for PMSM).

In Experiments~\ref{subsec:2d_burgers} and \ref{subsec:6d_burgers}, a Neumann boundary condition is additionally imposed during the training of the velocity network $\psi_{\boldsymbol{\omega}}$. Since the Burgers' equations feature infinitely long steep front structures, the samples may flow out of the computational domain $\Omega$ without boundary constraints, leading to a reduction of effective samples inside the domain at later times. To address this, we add a boundary penalty term to the loss function of the velocity potential network:
$$
\mathcal{J}_{\mathbf{V}}^{\text{bc}}({\boldsymbol{\omega}}) = \int_{\partial\Omega} \left( \frac{\partial \psi_{\boldsymbol{\omega}}}{\partial \mathbf{n}} \right)^2 dS = \int_{\partial\Omega} \left( \nabla \psi_{\boldsymbol{\omega}} \cdot \mathbf{n} \right)^2 dS ,
$$
where $\mathbf{n}$ is the outward unit normal vector on the boundary. This penalty term alleviates the issue of sampling points drifting out of the computational domain to some extent, allowing most of samples to concentrate within the interior.

To present the results in a clear and concise manner, we adopt a plotting interval of $0.1$, which differs from the computational time step $\Delta t = 0.05$. This means that not all adaptive sampling points at every time level are displayed, yet the selected results are sufficient to demonstrate the effectiveness of our method.

Table~\ref{tab:paper_result_summary} summarizes the best results of various methods across all experiments, highlighting the improvements achieved by our enhanced approach. Detailed descriptions of the experiments and comprehensive results are provided in the respective subsections.

\begin{table}[htbp]
\caption{Error comparison among PINNs, MSM, and PMSM across four benchmark problems.}
\centering
\small
\begin{tabular}{p{3.2cm}p{1.3cm}p{2.8cm}p{2.2cm}}
\toprule
Problem & Method & Relative $L^2$ error & $L^\infty$ error \\
\midrule
2D Burgers'  & PINNs & $1.359\times 10^{-3}$ & $3.611\times 10^{-2}$ \\
& MSM & $1.292\times 10^{-3}$ & $1.305\times 10^{-2}$ \\
& PMSM & $3.878\times 10^{-4}$ & $5.403\times 10^{-3}$ \\
\midrule
2D parabolic  & PINNs & $1.115\times 10^{-1}$ & $3.694\times 10^{-1}$ \\
& MSM & $8.118\times 10^{-2}$ &  $2.068\times 10^{-1}$ \\
& PMSM & $4.338\times 10^{-2}$ & $5.618\times 10^{-2}$ \\
\midrule
3D Fokker--Planck  & PINNs & $2.507\times 10^{-1}$ & $7.735\times 10^{-1}$ \\
& MSM & $1.495\times 10^{-2}$ & $1.157\times 10^{-1}$ \\
& PMSM & $7.267\times 10^{-3}$ & $3.127\times 10^{-2}$ \\
\midrule
6D Burgers'  & PINNs & $6.476\times 10^{-3}$ & $2.303\times 10^{-1}$ \\
& MSM &$7.243 \times 10^{-4}$ & $1.835\times 10^{-2}$ \\
& PMSM & $2.456\times 10^{-4}$ & $1.053\times 10^{-2}$ \\
\bottomrule
\end{tabular}
\label{tab:paper_result_summary}
\end{table}

\subsection{2D Burgers' equation} \label{subsec:2d_burgers}

Consider the following 2D Burgers' equation:
\begin{equation}\label{eq:Burgers}
\left\{
\begin{array}{ll}
    \frac{\partial u}{\partial t} = \alpha (\frac{\partial^2 u}{\partial x^2} + \frac{\partial^2 u}{\partial y^2}) - u(\frac{\partial u}{\partial x} + \frac{\partial u}{\partial y}), & \text{in}\ \Omega\times[0, T], \\
    u(x, y, 0) = \frac{1}{1+e^{\frac{x+y}{2\alpha}}},  & \text{in}\ \Omega,\\
    u(x, y, t) = \frac{1}{1+e^{\frac{x+y-t}{2\alpha}}}, & \text{on}\ \partial\Omega\times[0, T].
\end{array}\right.
\end{equation}
The exact solution is
\[
u(x, y, t) = \frac{1}{1+e^{\frac{x+y-t}{2\alpha}}}.
\]
In this experiment, we set $\alpha = 0.001$, $T = 1.05$, and $\Omega = [-1,1]^2$, so the front location is given by $x+y=t$.

This equation combines nonlinear convection with diffusion, whose solution develops sharper gradients and stronger nonlinear features as the diffusion parameter $\alpha$ decreases, thereby increasing the difficulty of approximation. Since the region with steep fronts is very small, the complex area is confined to a thin transition layer. Because of its limited size, uniform sampling hardly covers this region, leading to poor approximation of the solution around there.

PMSM and MSM use $N=500$ adaptive points and $N_u=1000$ uniform PDE points per time step, while the PINNs uses $1500$ fixed interior points. With $\Delta t =0.05$, we perform $20$ iterations in total, resulting in the final training time scale $\{0, 0.05, \dots, 1.05\}$. We set the network $u_{\bXheta}$ to be trained for $7500 + 20 \times 1500 + 15000 = 52500$ epochs in our method, where the three terms correspond to pretraining, extension-stage iterations, and final refinement, respectively. The PINNs and MSM baselines train the solution network for the same epochs.

Under the above parameter settings, the relative $L^2$ error decreases from $1.359\times 10^{-3}$ for PINNs to $3.878\times 10^{-4}$ for PMSM, and the $L^\infty$ error decreases from $3.611\times 10^{-2}$ to $5.403\times 10^{-3}$. The errors of MSM have improved compared to PINNs, but it is still higher than that of PMSM. This pattern is consistent with the geometry of the problem: the main benefit of PMSM is not a uniform improvement over the whole domain, but a more accurate and persistent resolution of the narrow traveling front, and the progressive extension strategy of PMSM amplifies this benefit.

\begin{figure}[htbp]
\centering
\includegraphics[width=0.88\textwidth]{\detokenize{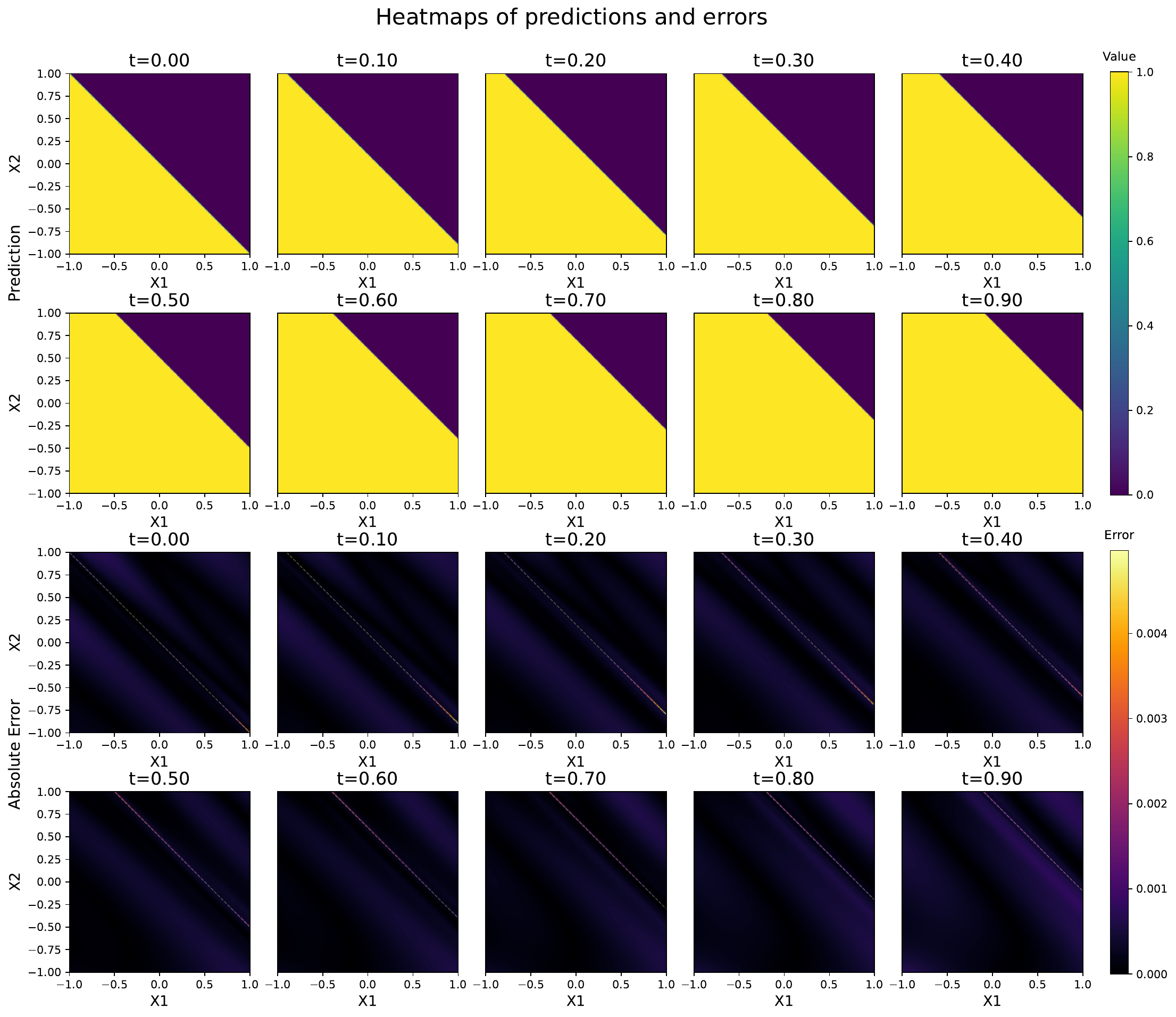}}
\caption{PMSM predicted solution and absolute error heatmaps for the 2D Burgers' equation \eqref{eq:Burgers}}
\label{fig:2d_burgers_heatmap}
\end{figure}

\begin{figure}[htbp]
\centering
\includegraphics[width=0.80\textwidth]{\detokenize{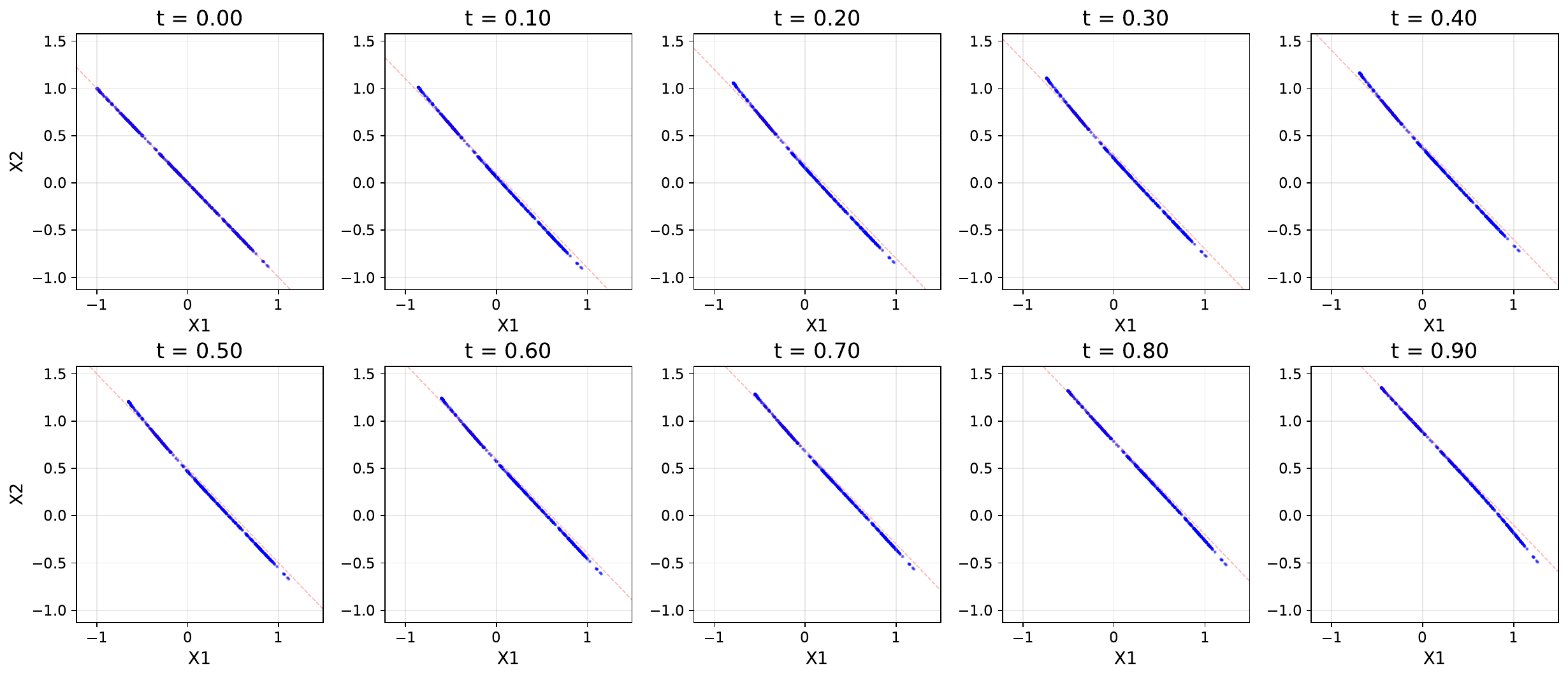}}
\caption{PMSM adaptive sample trajectories for the 2D Burgers' equation \eqref{eq:Burgers}}
\label{fig:2d_burgers_samples}
\end{figure}

Figure~\ref{fig:2d_burgers_heatmap} presents the results obtained by our method for this example together with the corresponding pointwise error, where one can observe that the error remains small even in regions with strong singularities. Figure~\ref{fig:2d_burgers_samples} shows the adaptive sampling points obtained by our method, which at all time concentrate around the theoretical wavefront $x+y=t$, ensuring accurate approximation in this region despite its narrow extent. These results demonstrate that our method achieves accurate approximation of the linear velocity field even under strongly singular conditions, further ensuring that the solution exhibits excellent accuracy.

\begin{figure}[htbp]
\centering
\includegraphics[width=0.88\textwidth]{\detokenize{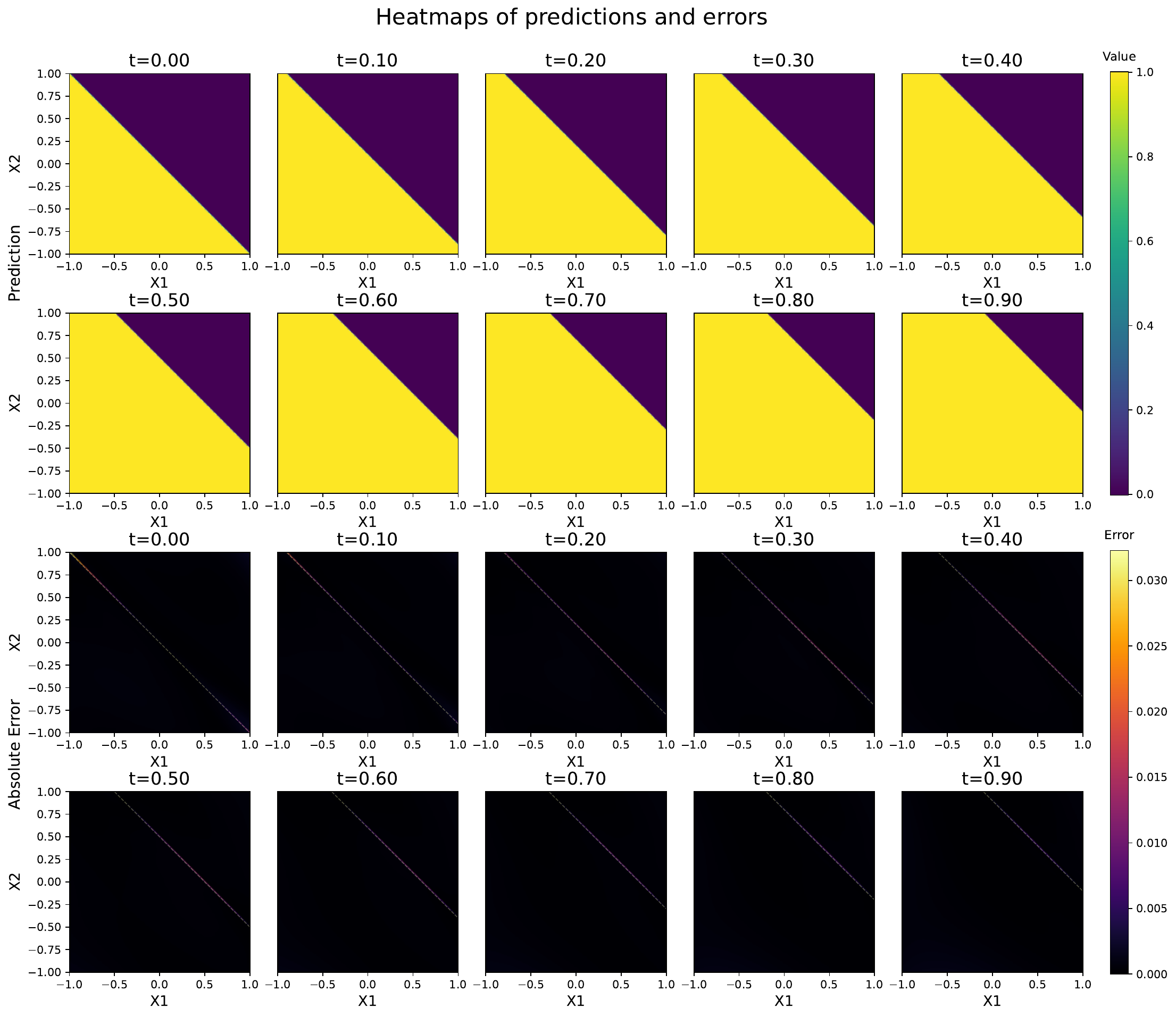}}
\caption{PINNs solution and absolute error heatmaps for the 2D Burgers' equation \eqref{eq:Burgers} }
\label{fig:2d_burgers_heatmap_PINN}
\end{figure}

\begin{figure}[htbp]
\centering
\includegraphics[width=0.88\textwidth]{\detokenize{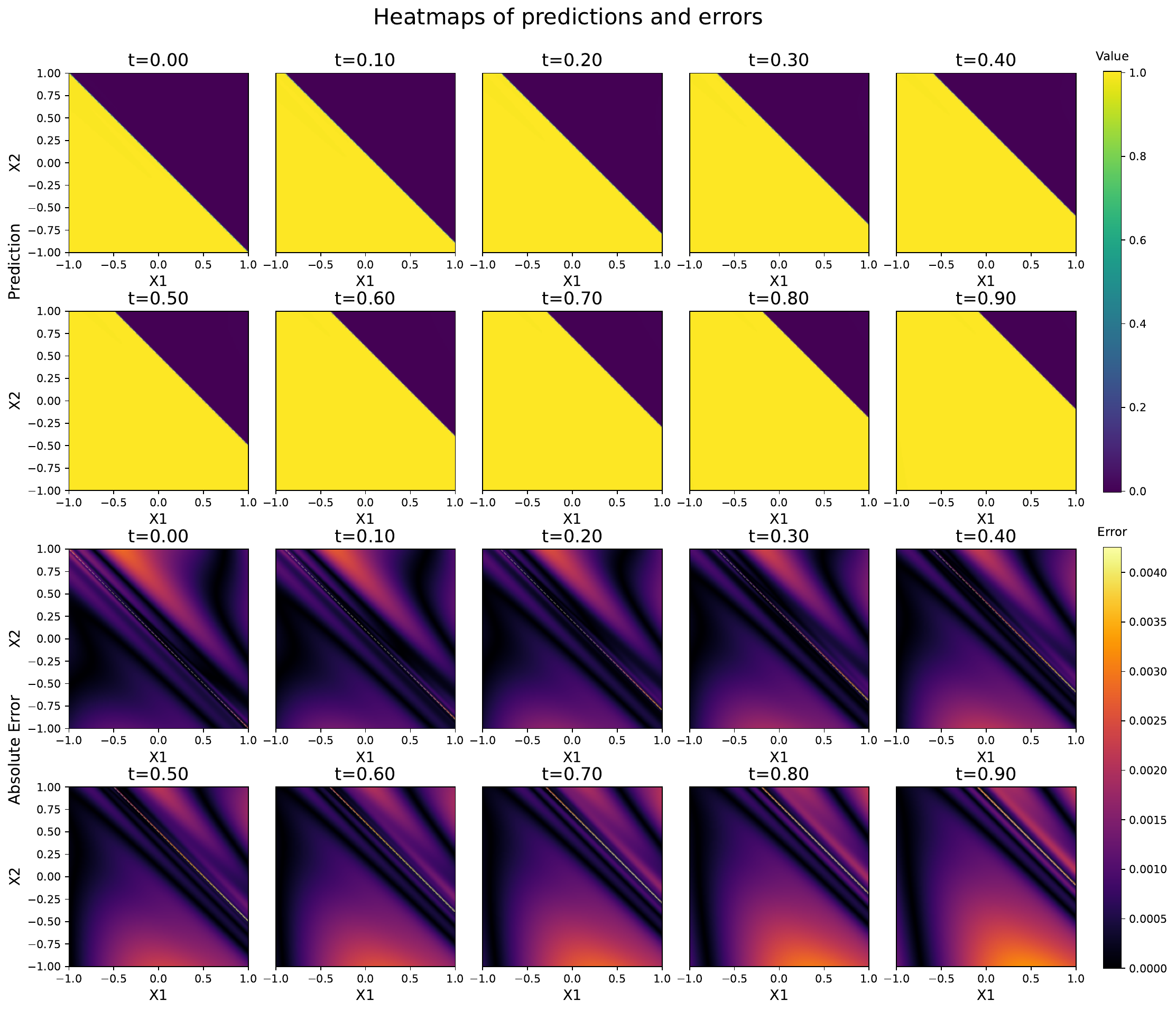}}
\caption{MSM solution and absolute error heatmaps for the 2D Burgers' equation \eqref{eq:Burgers}}
\label{fig:2d_burgers_heatmap_MSM}
\end{figure}

\begin{figure}[htbp]
\centering
\includegraphics[width=0.80\textwidth]{\detokenize{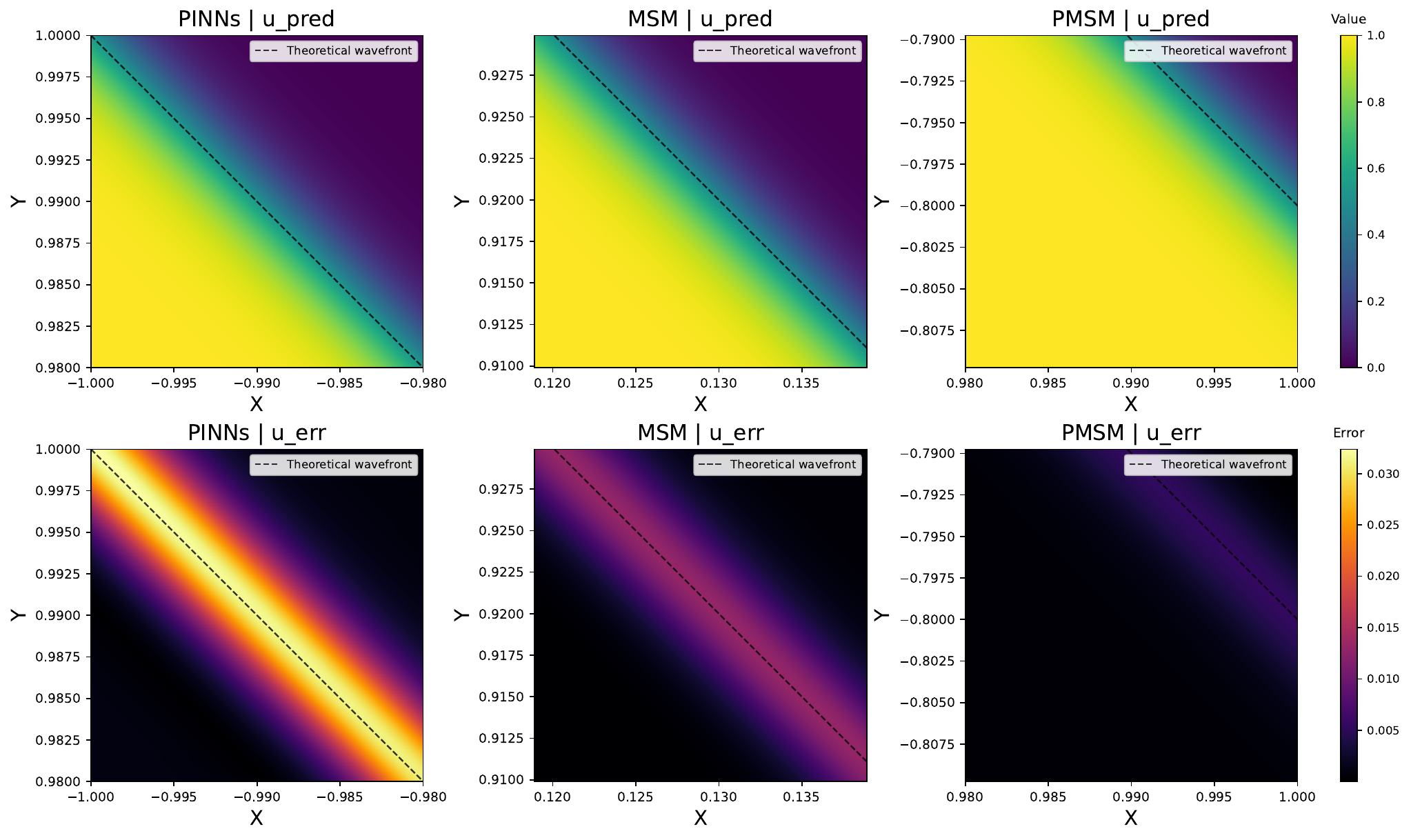}}
\caption{Enlarged views around the regions of maximum absolute error for the 2D Burgers' equation \eqref{eq:Burgers} obtained by PINNs (left), MSM (middle) and PMSM (right)}
\label{fig:2d_burgers_3contrast}
\end{figure}

Figure~\ref{fig:2d_burgers_heatmap_PINN} and \ref{fig:2d_burgers_heatmap_MSM} shows the results obtained by PINNs and MSM respectively. Although the global approximations look similar, the $L^\infty$ errors are significantly different, which means the location of the wavefront is not accurate for PINNs (see the left figure in Figure~\ref{fig:2d_burgers_3contrast}). Note that this result obtained by PINNs is the best selected from several runs, since we know the result. So uniformly sampled collocation points in the vanilla PINNs can not capture the key structure of the solution stably. In contrast, both MSM and PMSM yield more stable and accurate results, with PMSM exhibiting a better performance (see the right figure in Figure~\ref{fig:2d_burgers_3contrast}).

\subsection{2D parabolic equation}  \label{subsec:2d_parabolic}

Consider a 2D parabolic equation with a moving Gaussian peak whose amplitude grows and saturates:
\begin{equation}\label{eq:parabolic4}
    \left\{
    \begin{array}{ll}
        \frac{\partial u}{\partial t} - \frac{\partial ^2u}{\partial x^2} - \frac{\partial ^2u}{\partial y^2} = f(x, y, t), & \text{in}\ \Omega\times[0, T],\\
        u(x, y, t) = u_b(x, y, t), & \text{on}\ \partial\Omega\times [0, T], \\
        u(x, y, 0) = C_0 e^{-\frac{x^2+y^2}{\alpha}}, & \text{in}\ \Omega .
    \end{array}
    \right.
\end{equation}
The exact solution is
\[
u(x, y, t) = \Big(C_{\infty}-(C_{\infty}-C_0)e^{-\beta t}\Big)e^{- \frac{(x-t)^2 + (y-t^2)^2}{\alpha}},
\qquad \beta>0.
\]
We use $\Omega=[-0.2,2.5]^2$, $T = 1.55$, $\alpha=0.01$, $C_{\infty}=4.0$, $C_0=1.0$, and $\beta=2.0$. The peak center moves along the curved trajectory $(x, y) = (t,t^2)$, its direction changes continuously with time, and its amplitude also varies. Compared with the Burgers' equations, this problem is less about following a thin front and more about keeping the complex region covered while its position, orientation, and magnitude evolve together. This is a different kind of examination for sample motion: points that are well placed at an early time still need to bend with the trajectory and remain useful when the peak becomes larger at later times.

For this experiment, PMSM and MSM use $N=500$ adaptive points and $N_u=500$ uniform PDE points per time step, while the PINNs uses $1000$ fixed interior points. The rollout covers $t=0, 0.05, \ldots, 1.55$, corresponding to $30$ extension rounds. Since this experiment uses $30$ extension rounds, the total epochs is $7500 + 30\times1500 + 15000 = 67500$ for $u_{\bXheta}$ in PMSM. The PINNs and MSM baselines are trained for the same total number of solution network epochs. On this problem, MSM yields only a marginal improvement over PINNs: the relative $L^2$ error stays at $8.118\times 10^{-2}$ (versus $1.115\times 10^{-1}$ for PINNs), and the $L^\infty$ error drops from $3.694\times 10^{-1}$ to $2.068\times 10^{-1}$. This limited gain reflects the difficulty of the full-time-domain iterative strategy when the peak simultaneously changes its position, direction, and amplitude. Re-evolving samples from $t=0$ at each iteration may not adapt quickly enough to the combined variation. By contrast, PMSM reduces the relative $L^2$ error to $4.338\times 10^{-2}$ and the $L^\infty$ error to $5.618\times 10^{-2}$. The progressive time extension strategy, which advances one time step at a time and updates the velocity field using local residual information, is better able to follow the curved trajectory and respond to the growing peak amplitude.

\begin{figure}[htbp]
\centering
\includegraphics[width=0.88\textwidth]{\detokenize{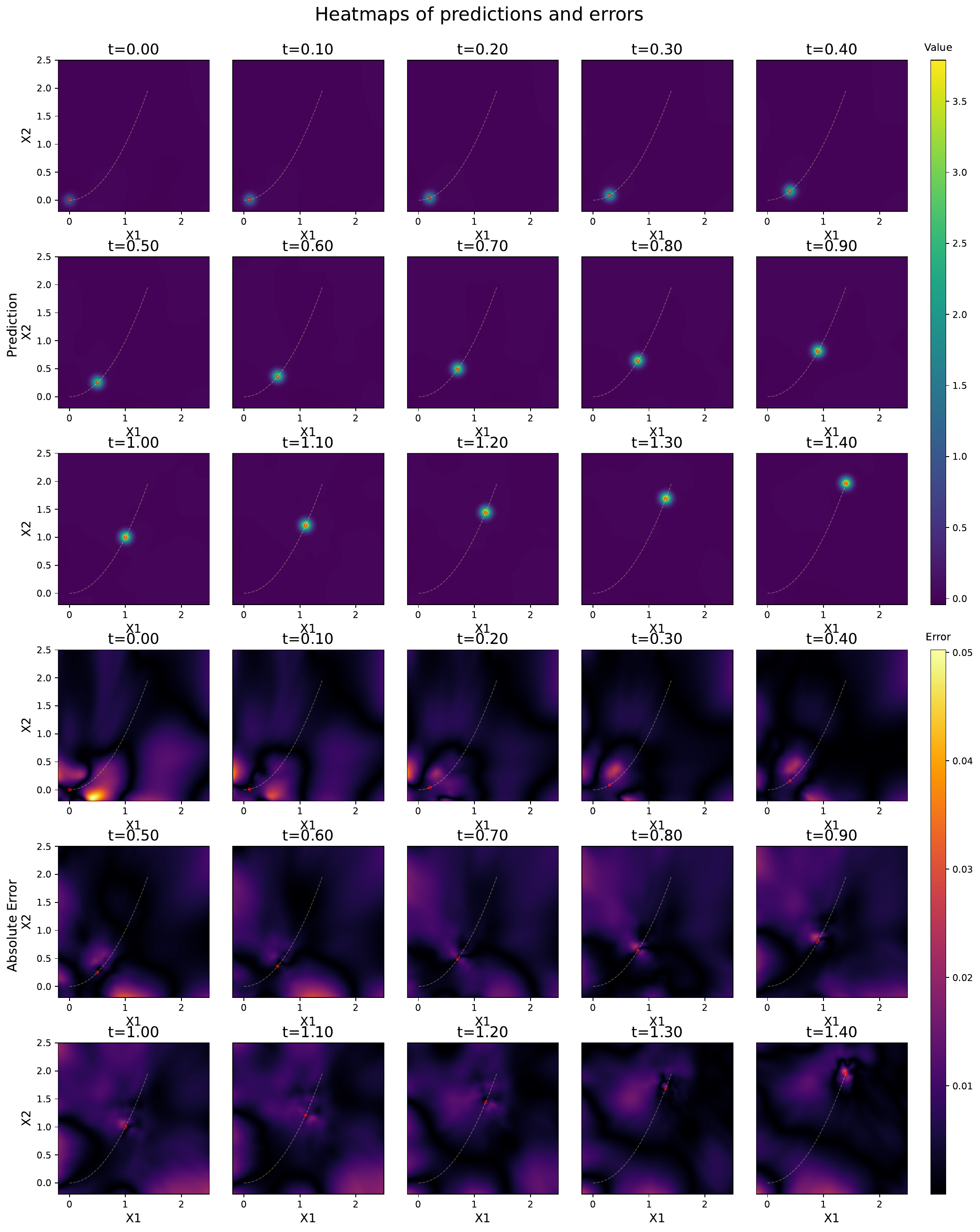}}
\caption{PMSM predicted solution (top) and absolute error (bottom) heatmaps for the 2D parabolic equation \eqref{eq:parabolic4} with a moving and amplifying Gaussian peak}
\label{fig:parabolic_heatmap}
\end{figure}

\begin{figure}[htbp]
\centering
\includegraphics[width=0.80\textwidth]{\detokenize{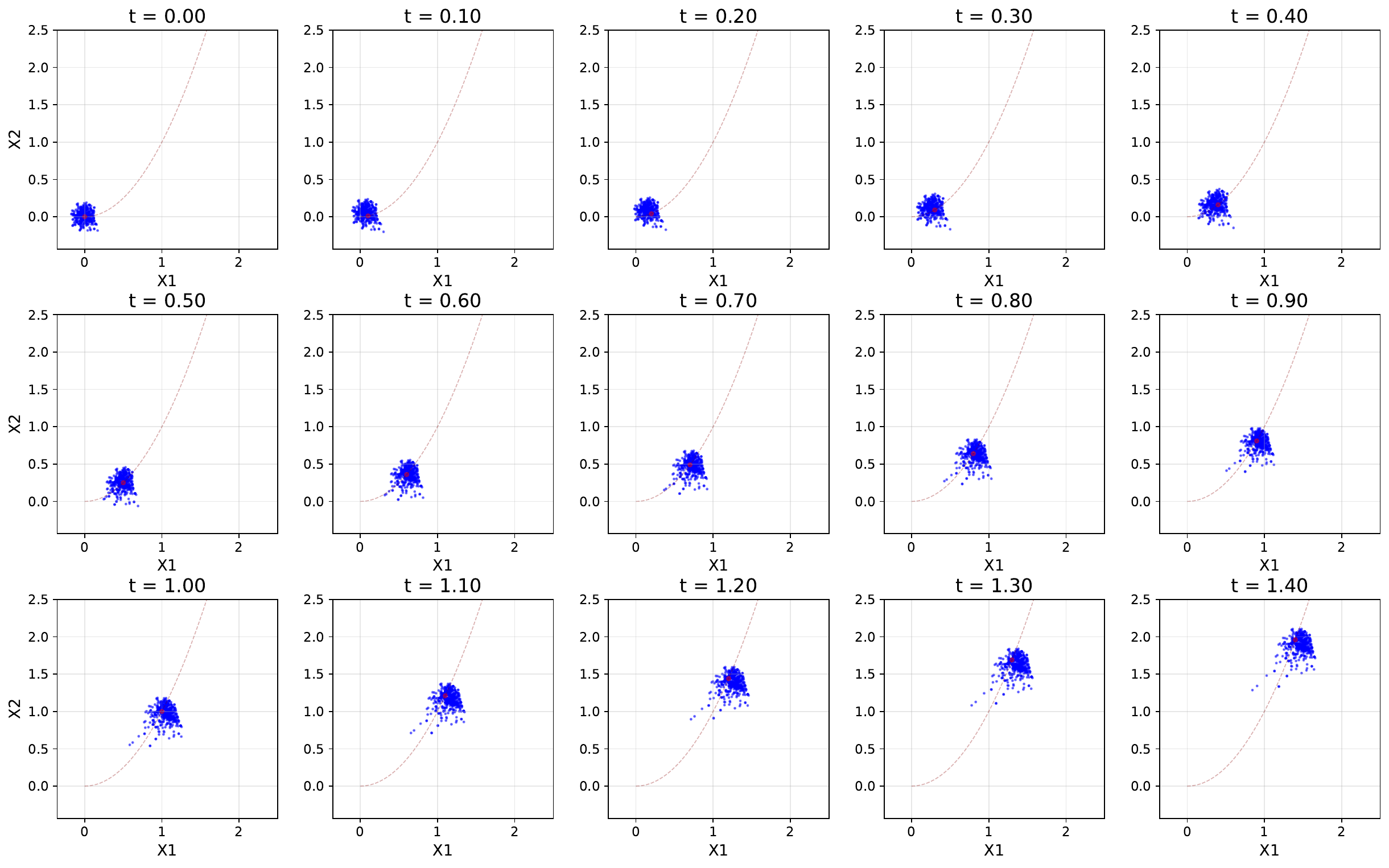}}
\caption{PMSM adaptive sample trajectories for the 2D parabolic equation \eqref{eq:parabolic4}}
\label{fig:parabolic_samples}
\end{figure}

\begin{figure}[htbp]
\centering
\includegraphics[width=0.88\textwidth]{\detokenize{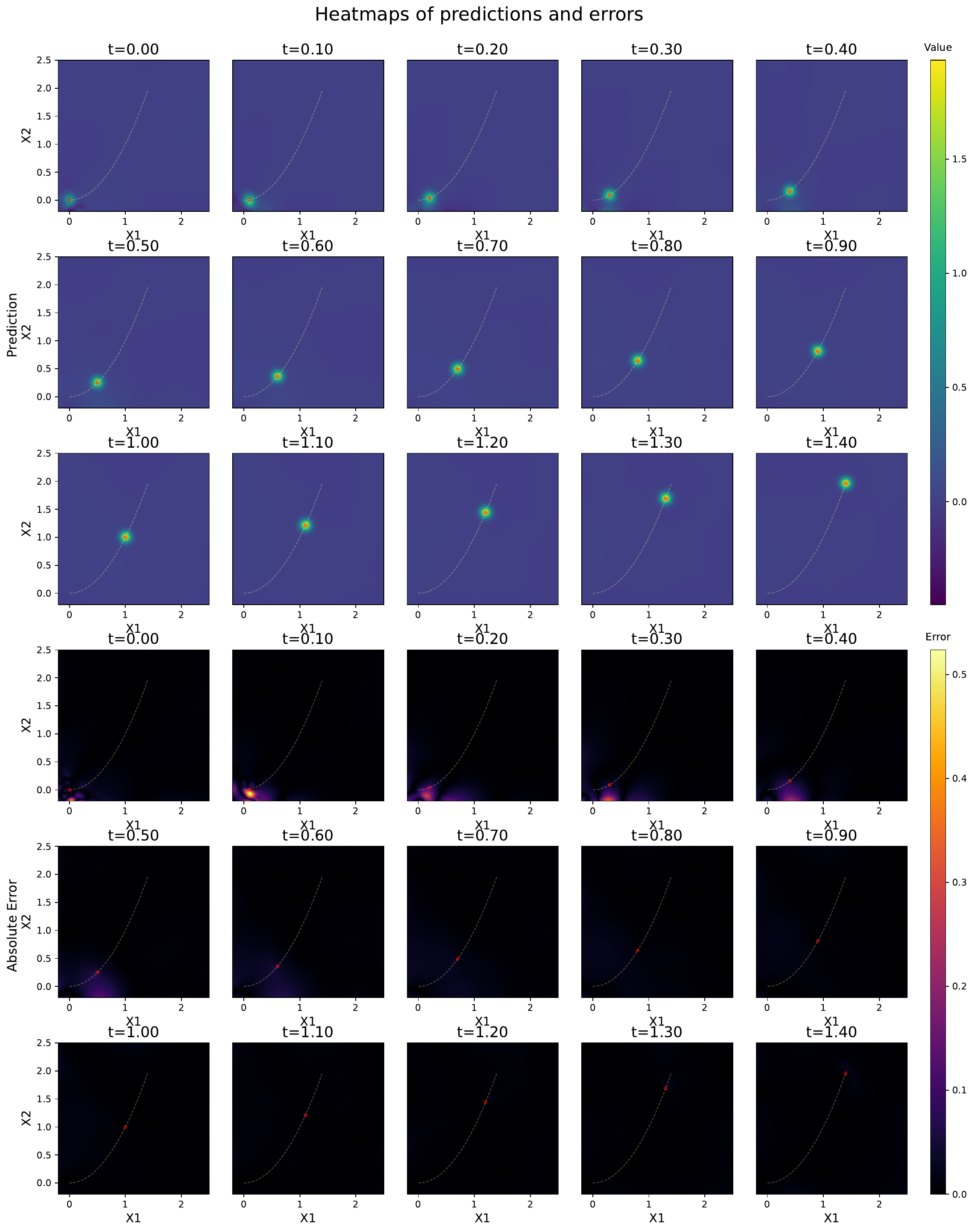}}
\caption{PINNs solution (top) and absolute error (bottom) heatmaps for the 2D parabolic equation \eqref{eq:parabolic4}}
\label{fig:parabolic_heatmap_PINN}
\end{figure}

\begin{figure}[htbp]
\centering
\includegraphics[width=0.88\textwidth]{\detokenize{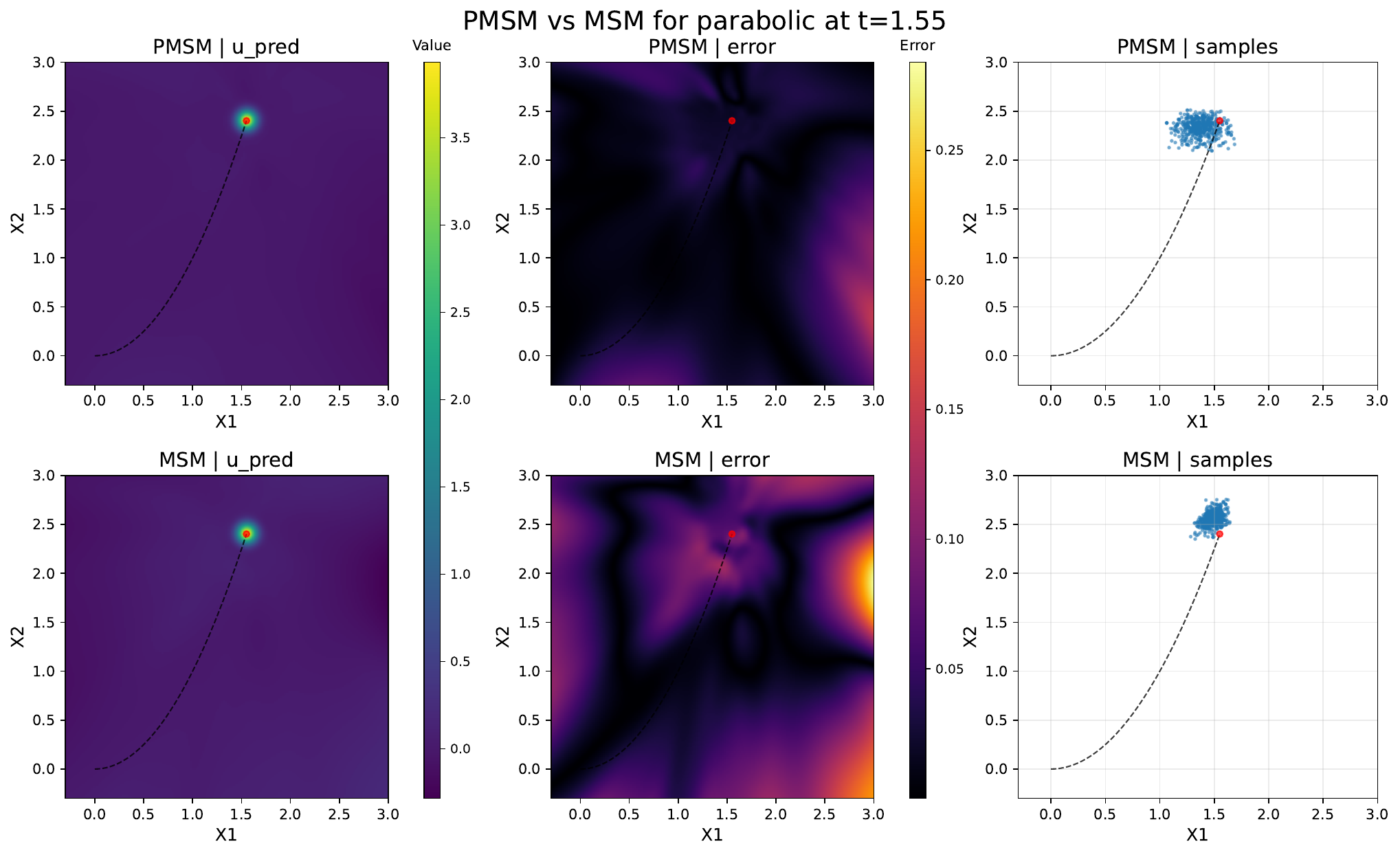}}
\caption{
PMSM (top) and MSM (bottom) solution heatmaps (left), absolute error heatmaps (middle) and adaptive sample trajectories (right) at $t=1.55$ for the 2D parabolic equation \eqref{eq:parabolic4}
}
\label{fig:parabolic_vs_t155}
\end{figure}

Figure~\ref{fig:parabolic_heatmap} shows that the dominant error is concentrated near the peak trajectory and in later times where the amplitude becomes larger. Figure~\ref{fig:parabolic_samples} shows that PMSM keeps a high-density sample cloud near the moving peak and gradually reshapes that cloud as the trajectory bends and the amplitude grows. For comparison, Figure~\ref{fig:parabolic_heatmap_PINN} shows the results obtained by PINNs. For PINNs method, the absolute error is noticeably large (see the bottom figure in Figure~\ref{fig:parabolic_heatmap_PINN}), indicating that it lacks sufficient accuracy in fitting the example with nonlinear moving and growing singularity. In Figure~\ref{fig:parabolic_vs_t155}, we present a comparison at $t=1.55$, showing the heatmaps of the solutions obtained by PMSM and MSM, together with their respective errors and the corresponding adaptive sampling points. Since MSM computes the velocity field over the entire time interval, numerical errors may accumulate when evolving the samples with this fixed velocity field. Moreover, the training of a potential field network with time-dependent magnitude and nontrivial variation inherently poses further challenges. As a result, the sampling points fail to move along with the singularity for large $t$, which in turn affects the overall solution quality and stability to a certain extent. However, PMSM performs better and still keeps a certain number of collocation points in the singularity region. This example therefore highlights that PMSM performs better for singularities with more complicated time dependence.

\subsection{3D Fokker-Planck equation}
\label{subsec:3d-fokker-planck}

Consider the following 3D SDE
\begin{equation}\label{eq:SDE}
    \left\{
    \begin{aligned}
        dx &= (-4(x-e^{-t})((x-e^{-t})^2 + (y-e^{-t})^2 + (z-e^{-t})^2 - r^2) - e^{-t})dt + \sigma dW_1, \\
        dy &= (-4(y-e^{-t})((x-e^{-t})^2 + (y-e^{-t})^2 + (z-e^{-t})^2 - r^2) - e^{-t})dt + \sigma dW_2, \\
        dz &= (-4(z-e^{-t})((x-e^{-t})^2 + (y-e^{-t})^2 + (z-e^{-t})^2 - r^2) - e^{-t})dt + \sigma dW_3, \\
    \end{aligned}
    \right.
\end{equation}
where $W_1, W_2, W_3$ denote independent Wiener processes. The corresponding Fokker-Planck equation is
\begin{equation}\label{eq:FP}
        \frac{\partial u}{\partial t} = \mathcal{L}u = -\nabla \cdot (fu) + \nabla^2 (Du),
\end{equation}
where $f = \begin{pmatrix}
    -4(x-e^{-t})((x-e^{-t})^2 + (y-e^{-t})^2 + (z-e^{-t})^2 - r^2) - e^{-t} \\
    -4(y-e^{-t})((x-e^{-t})^2 + (y-e^{-t})^2 + (z-e^{-t})^2 - r^2) - e^{-t} \\
    -4(z-e^{-t})((x-e^{-t})^2 + (y-e^{-t})^2 + (z-e^{-t})^2 - r^2) - e^{-t} \\
\end{pmatrix}, D=\frac{\sigma^2}{2}\mathbf{I}$. The initial condition is
\[
u(x, y, z, 0) = \frac{1}{K} e^{-\frac{2}{\sigma^2}((x-1)^2 + (y-1)^2 + (z-1)^2 - r^2)^2},
\]
where
\[
K=\iiint e^{-\frac{2}{\sigma^2}((x-1)^2 + (y-1)^2 + (z-1)^2 - r^2)^2}\,dxdydz.
\]
The exact solution is therefore
\[
u(x, y, z, t) = \frac{1}{K} e^{-\frac{2}{\sigma^2}((x-e^{-t})^2 + (y-e^{-t})^2 + (z-e^{-t})^2 - r^2)^2},
\]
and the Dirichlet boundary condition is imposed from this analytical solution.

This equation describes the evolution of a probability density under a 3D Fokker-Planck operator, where the drift field $f$ induces a nonlinear transport toward a moving spherical structure. As time evolves, the shell center moves from $(1,1,1)$ toward the origin along the diagonal direction, while the density remains strongly localized around the radius-$r$ surface. The main numerical difficulty is therefore not only the three-dimensionality of the domain, but also the fact that the informative region occupies a small fraction of the volume and has a nontrivial hollow geometry. Uniform collocation tends to waste many points in low density interior and exterior regions, and boundary data alone do not resolve where the probability mass should sit inside the domain. An adaptive method should instead keep concentrating samples near the moving high-probability shell throughout the evolution.

In this experiment, we use $\sigma=0.1$, $r=1.0$, $\Omega=[0,2.1]^3$ and $T = 1.15$. With $\Delta t=0.05$, the rollout starts from $6$ initial time steps and appends one new time step at each of the $18$ extension rounds. Under this rollout, the PMSM solution network is trained for $30000 + 18\times6000 + 60000 = 198000$ epochs. The PINNs and MSM baselines are trained for the same $198000$ epochs. PMSM uses $N=1000$ adaptive points and $N_u=2500$ uniform PDE points per time step, together with $N_0=1000$ adaptive initial-condition points, MSM uses the same as PMSM, while PINNs uses $3500$ fixed interior points and the same $1000$ initial points.

As in the previous examples, MSM yields errors between those of PINNs and PMSM. PMSM reduces the relative $L^2$ error to $7.267\times 10^{-3}$ and the $L^\infty$ error to $3.127\times 10^{-2}$, compared with $2.507\times 10^{-1}$ and $7.735\times 10^{-1}$ for PINNs. The gain is meaningful in view of the shell geometry: a large $L^\infty$ error typically indicates that the learned density fails to place mass accurately near the narrow concentration surface, while the smaller PMSM error suggests that the moving samples are able to follow that surface more reliably.

\begin{figure}[htbp]
\centering
\includegraphics[width=0.88\textwidth]{\detokenize{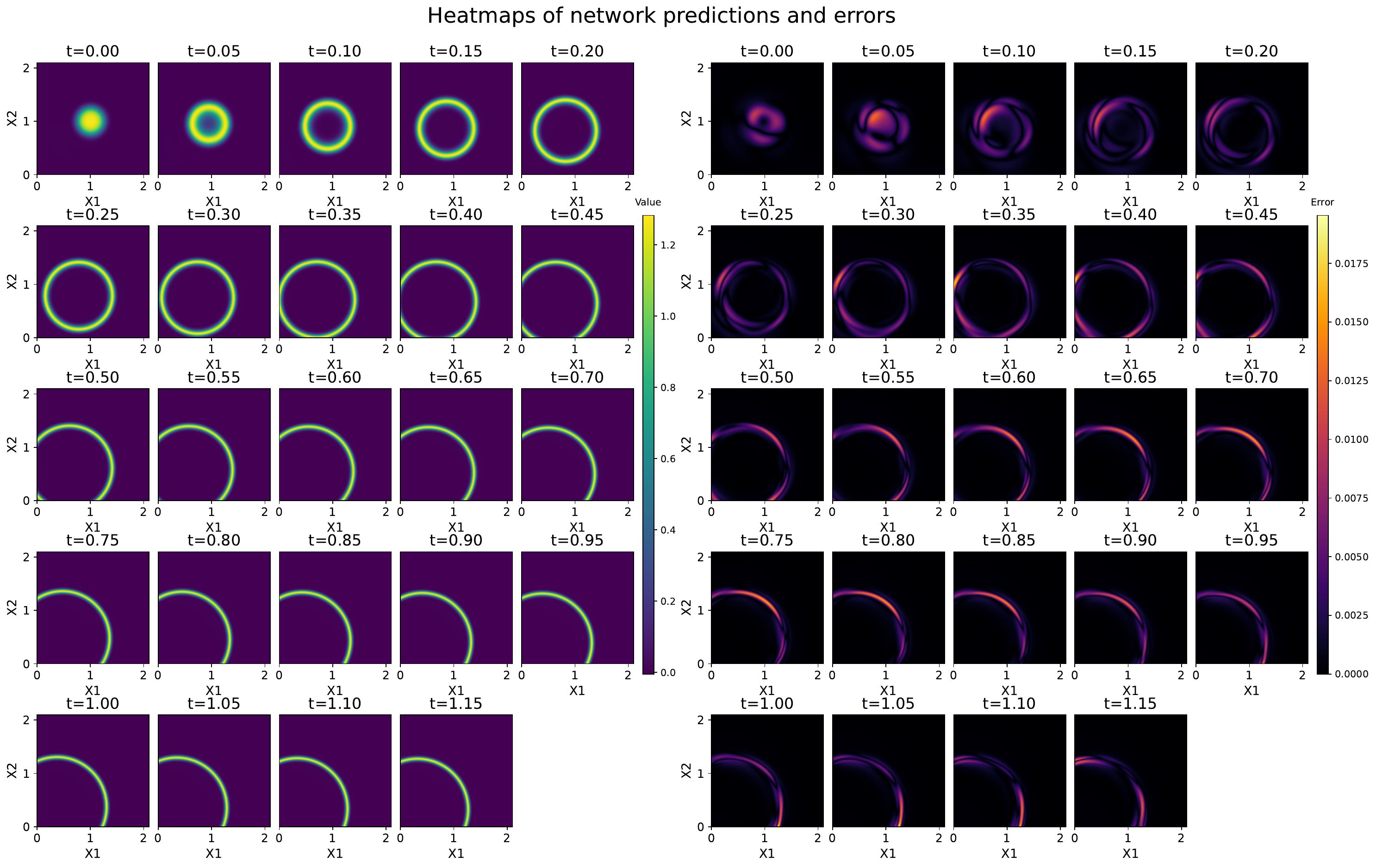}}
\caption{PMSM two-dimensional heatmap of the solution on the slice $z=0$ (left) and absolute error (right) visualizations for the 3D Fokker--Planck equation \eqref{eq:FP}}
\label{fig:fokkerplanck_heatmap}
\end{figure}

\begin{figure}[htbp]
\centering
\includegraphics[width=0.80\textwidth]{\detokenize{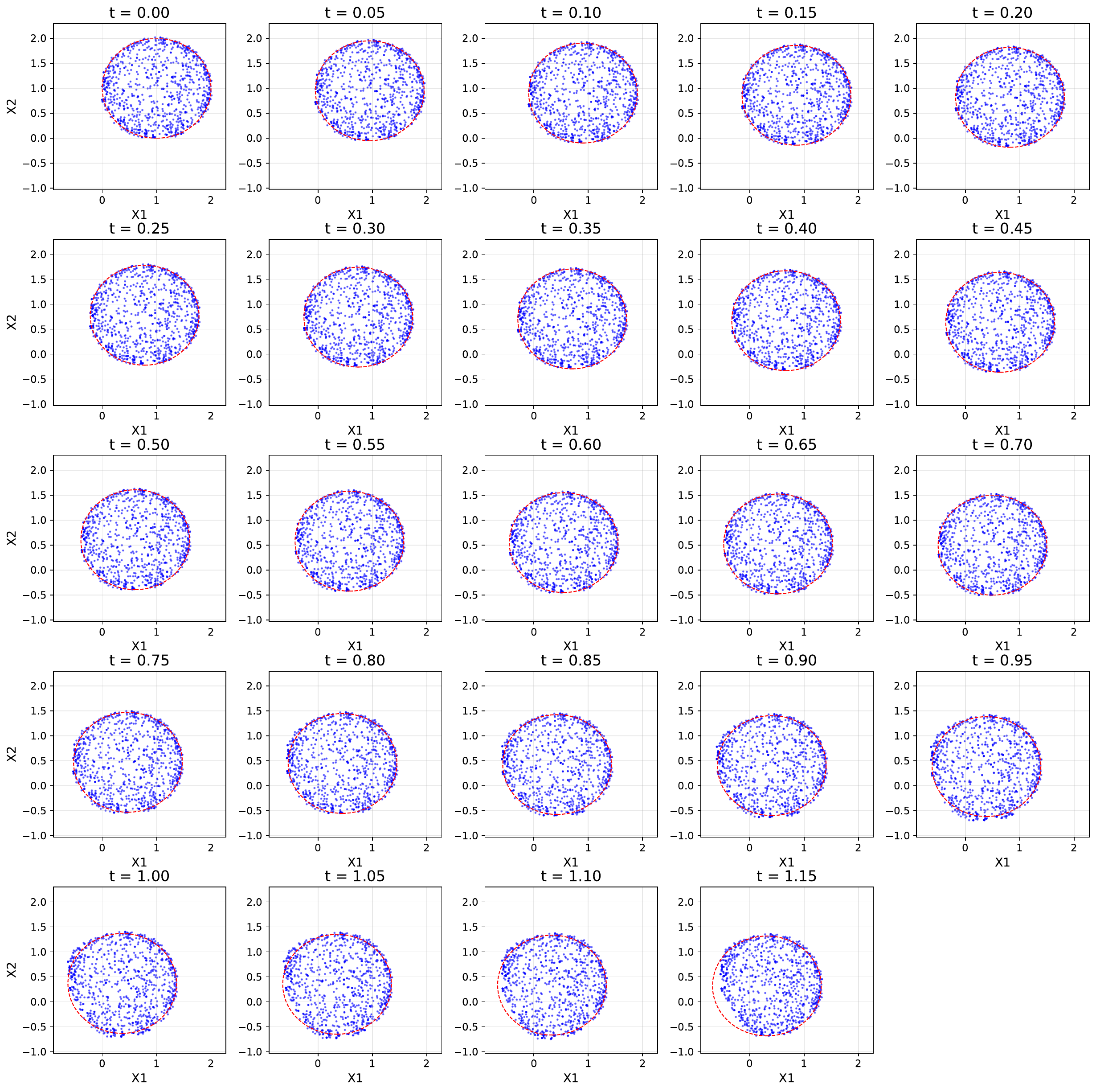}}
\caption{Projection of adaptive sampling trajectories for the 3D Fokker--Planck equation \eqref{eq:FP} generated by PMSM in two dimensions. For better visualization, the plotting domain is enlarged compared with the computational domain}
\label{fig:fokkerplanck_samples}
\end{figure}

\begin{figure}[htbp]
\centering
\includegraphics[width=0.88\textwidth]{\detokenize{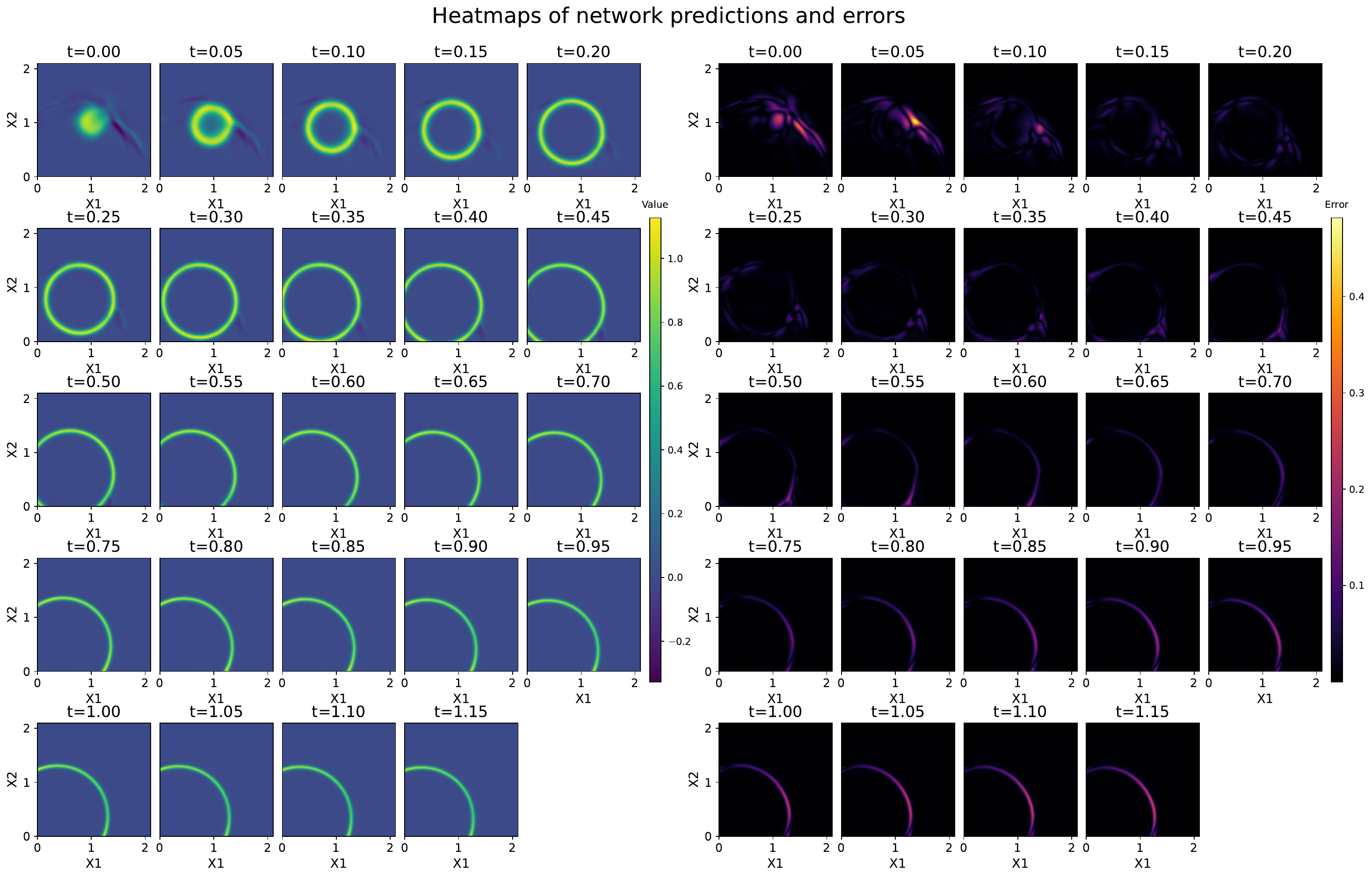}}
\caption{PINNs two-dimensional heatmap of the 3D solution on the slice $z=0$ (left) and absolute error (right) visualizations for the 3D Fokker--Planck equation \eqref{eq:FP}}
\label{fig:fokkerplanck_heatmap_PINN}
\end{figure}

\begin{figure}[htbp]
\centering
\includegraphics[width=0.88\textwidth]{\detokenize{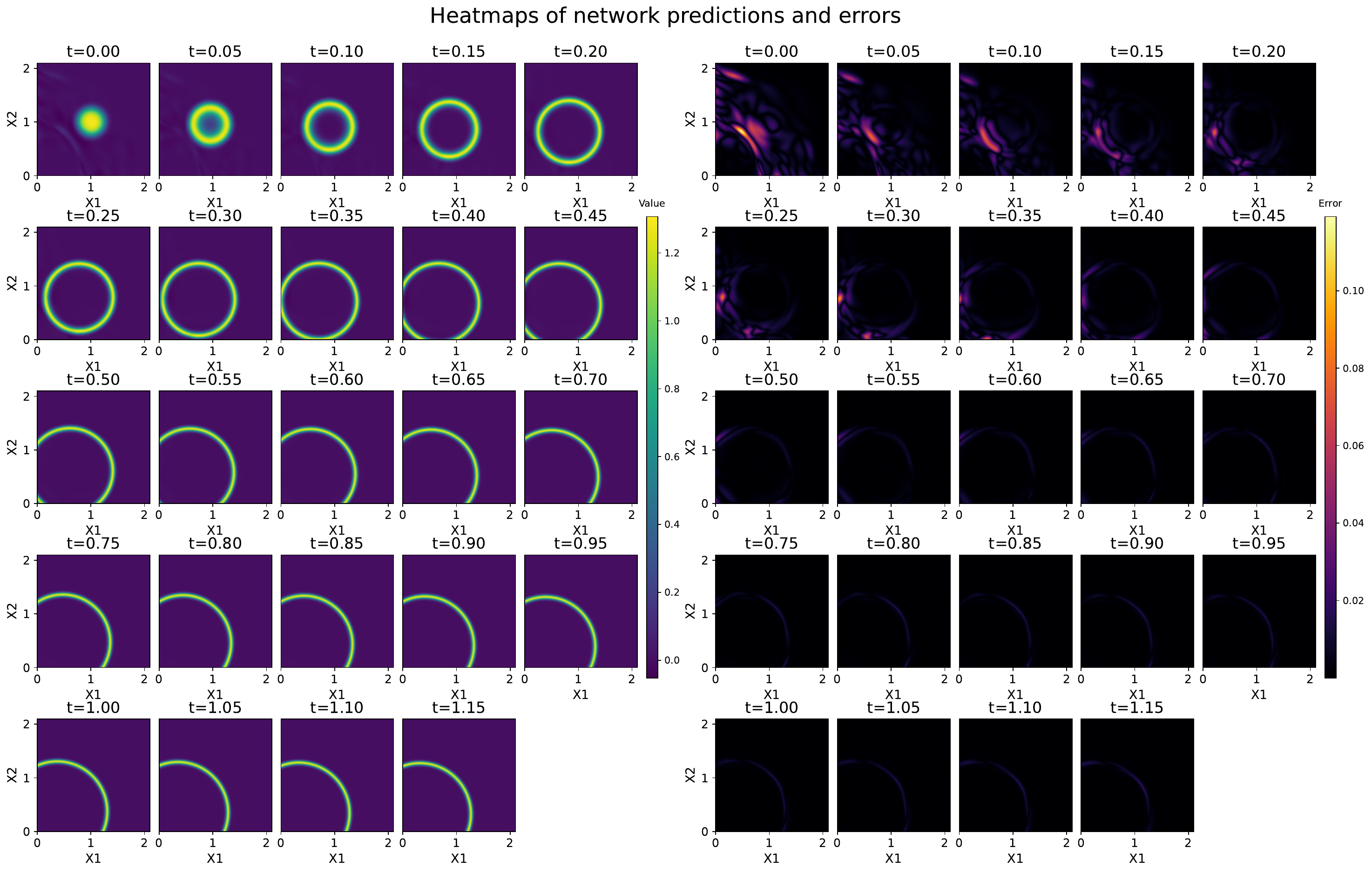}}
\caption{MSM two-dimensional heatmap of the 3D solution on the slice $z=0$ (left) and absolute error (right) visualizations for the 3D Fokker--Planck equation \eqref{eq:FP}}
\label{fig:fokkerplanck_heatmap_MSM}
\end{figure}

\begin{figure}[htbp]
\centering
\includegraphics[width=0.80\textwidth]{\detokenize{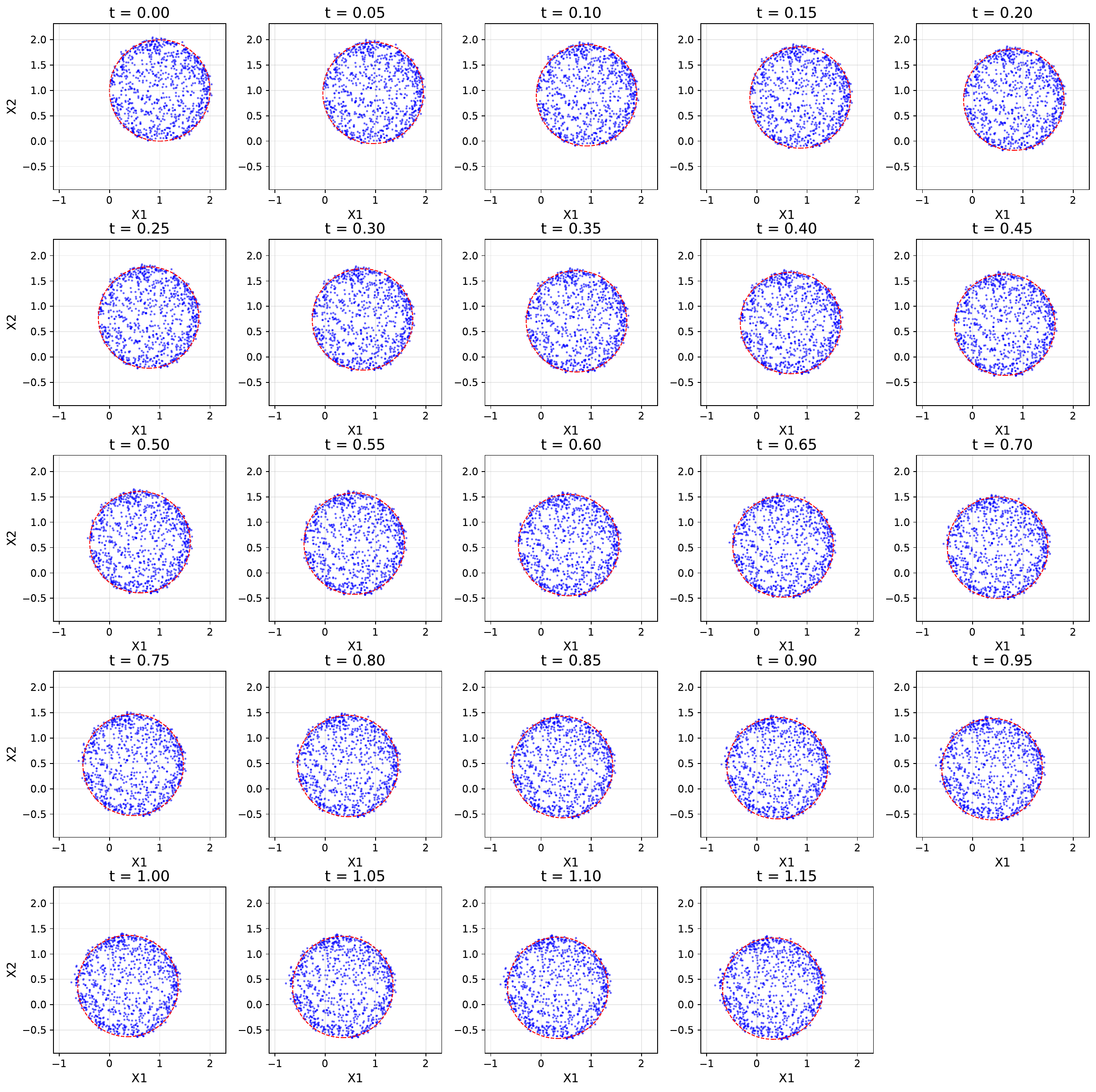}}
\caption{Projection of adaptive sampling trajectories for the 3D Fokker--Planck equation \eqref{eq:FP} generated by MSM in two dimensions. For better visualization, the plotting domain is enlarged compared with the computational domain}
\label{fig:fokkerplanck_samples_MSM}
\end{figure}

\begin{figure}[htbp]
\centering
\includegraphics[width=0.80\textwidth]{\detokenize{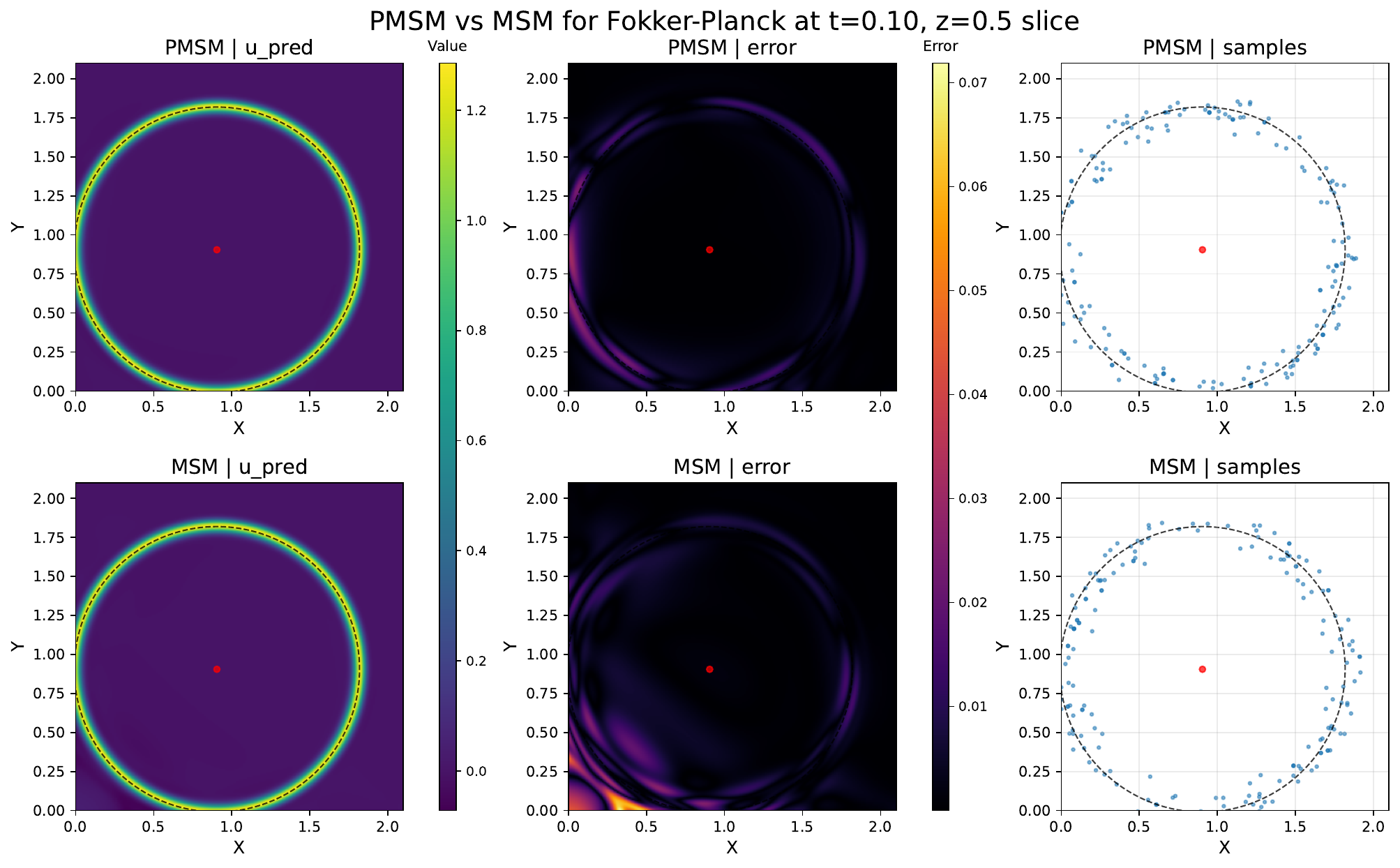}}
\caption{
PMSM (top) and MSM (bottom) solution heatmaps (left), absolute error heatmaps (middle) and projection of adaptive sample trajectories (right) at $t=0.1$, $z=0.5$ slice for the 3D Fokker--Planck equation \eqref{eq:FP}
}
\label{fig:fokkerplanck_vs_t010}
\end{figure}

Figure~\ref{fig:fokkerplanck_heatmap} shows that the dominant approximation error is concentrated near the narrow moving shell, while the rest of the domain contributes comparatively little to the overall difficulty. Figure~\ref{fig:fokkerplanck_samples} shows that PMSM responds by transporting the adaptive points together with the shell-shaped high density region instead of leaving large numbers of points in the low density interior and exterior. This behavior is consistent with the strong reduction in both relative $L^2$ and $L^\infty$ errors, and it further illustrates that PMSM can adapt not only to fronts and peaks, but also to more intricate moving geometric structures.

The results of PINNs and MSM are shown in Figure~\ref{fig:fokkerplanck_heatmap_PINN} and \ref{fig:fokkerplanck_heatmap_MSM}, while the adaptive samples obtained by MSM are given in Figure~\ref{fig:fokkerplanck_samples_MSM}. PINNs cannot accurately fit the distribution over the entire spherical surface, showing significant errors both at the initial time and for larger $t$, and the quantitative error for PINNs is significantly large, with the maximum absolute error reaching $7.735\times 10^{-1}$ (see the right figure in Figure~\ref{fig:fokkerplanck_heatmap_PINN}). In comparison, MSM yields a clear improvement: its sample trajectories are able to follow the moving shell region (see Figure~\ref{fig:fokkerplanck_samples_MSM}). However, because the solution and velocity field are learned across the entire time domain and the equation undergoes rapid changes during the initial period, the MSM error becomes sensitive to the weighting of the loss function and is still noticeable in the first several time steps. PMSM further reduces this early error through the progressive time-stepping strategy, which resolves the initial transient more carefully before extending to later times. To be specific, in Figure~\ref{fig:fokkerplanck_vs_t010}, we present a comparison at $t=0.1, z=0.5$, showing the heatmaps of the solutions, errors and the projection of corresponding adaptive sampling points obtained by PMSM and MSM, where the projection is taken over samples with $z \in [0.3, 0.7]$. In contrast, even with similarly distributed sample points, our method begins training over shorter time intervals and gradually extends the temporal domain. This strategy alleviates the difficulties of training at early time and leads to a noticeable improvement in accuracy.

\subsection{Higher-Dimensional Burgers' equation} \label{subsec:6d_burgers}

We next consider a variable-speed $d$-dimensional Burgers' equation and instantiate it in six dimensions:
\begin{equation}\label{eq:dBurgers_general}
\left\{
\begin{array}{ll}
\displaystyle
\frac{\partial u}{\partial t}
=
g'(t)\left(\alpha \sum_{i=1}^d \frac{\partial^2 u}{\partial x_i^2}
- u \sum_{i=1}^d \frac{\partial u}{\partial x_i}\right),
& \text{in}\ \Omega\times[0, T], \\[1.2ex]
\displaystyle
u(\boldsymbol{x}, 0)
=
\frac{1}{1+\exp\left(\frac{\sum_{i=1}^d x_i - \frac{d}{2}g(0)}{2\alpha}\right)},
& \text{in}\ \Omega, \\[1.2ex]
\displaystyle
u(\boldsymbol{x}, t)
=
\frac{1}{1+\exp\left(\frac{\sum_{i=1}^d x_i-\frac{d}{2}g(t)}{2\alpha}\right)},
& \text{on}\ \partial\Omega\times[0, T],
\end{array}
\right.
\end{equation}
where $g$ is continuously differentiable and monotonically increasing with respect to $t$. Its exact solution is
\[
u(\boldsymbol{x}, t)
=
\frac{1}{1+\exp\left(
\frac{\sum_{i=1}^d x_i-\frac{d}{2}g(t)}{2\alpha}
\right)}.
\]
In this experiment, we use $d=6$, $\Omega=[-3,3]^6$, $T=1.05$, $\alpha=0.01$, and the rational deceleration law
\[
g(t)=\frac{t}{1+ t},
\qquad
g'(t)=\frac{1}{(1+ t)^2}.
\]
The front therefore satisfies $\sum_{i=1}^6 x_i = 3g(t)$ and moves quickly at early times before gradually slowing down.

This equation represents a high-dimensional generalization of the Burgers equation, where the dynamics are governed simultaneously by nonlinear convection and diffusion across $d$ spatial directions. In six dimensions, the front is a high-dimensional transition surface, yet its neighborhood occupies an even smaller fraction of the ambient volume. The singularity region is therefore sparse from the viewpoint of random collocation, and a large portion of fixed points can fall far away from the part of the solution that determines the error. The rational law $g(t)$ adds a second difficulty: the theoretical velocity field is no longer uniform but exhibits nonlinear variations, so a useful sampler must learn the pace of the motion, not only its direction. PMSM addresses both aspects by retaining adaptive trajectories near the front and updating their motion through the residual time evolution.

We use $N=500$ adaptive points and $N_u=1000$ uniform PDE points for PMSM and MSM, while $1500$ fixed interior points for the PINNs. We follow the same training setup as in the 2D Burgers experiment \ref{subsec:2d_burgers}.

Under this high-dimensional and nonuniform setting, the relative $L^2$ error decreases from $6.476\times 10^{-3}$ for PINNs to $2.456\times 10^{-4}$ for PMSM. The $L^\infty$ error decreases from $2.303\times 10^{-1}$ to $1.053\times 10^{-2}$. The errors of MSM are like the previous example, between that of PINNs and PMSM. These are consistent with the geometric difficulty of the problem: the moving front occupies a small part of the six-dimensional domain, and transported samples provide denser coverage of that region. Since our PMSM does not train the velocity field over the entire time interval $[0, T]$ like MSM, but instead computes the field between two adjacent time at each extension step, the training becomes more efficient and accurate.

\begin{figure}[htbp]
\centering
\includegraphics[width=0.88\textwidth]{\detokenize{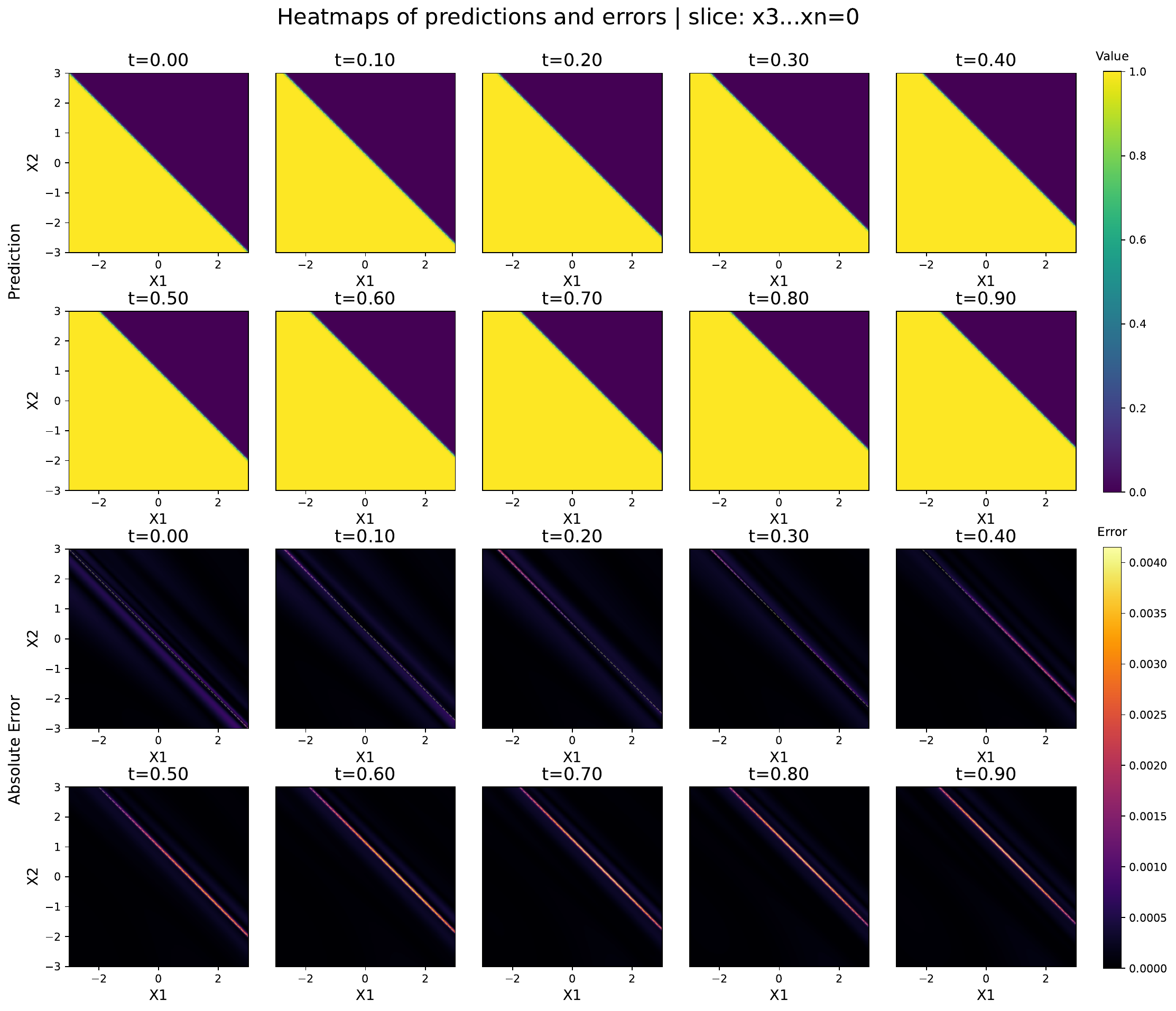}}
\caption{PMSM predicted solution (top) and absolute error (bottom) slice heatmaps for the 6D variable-speed Burgers' equation \eqref{eq:dBurgers_general}. The visualization uses the slice $x_3=x_4=x_5=x_6=0$}
\label{fig:6d_burgers_heatmap}
\end{figure}

\begin{figure}[htbp]
\centering
\includegraphics[width=0.80\textwidth]{\detokenize{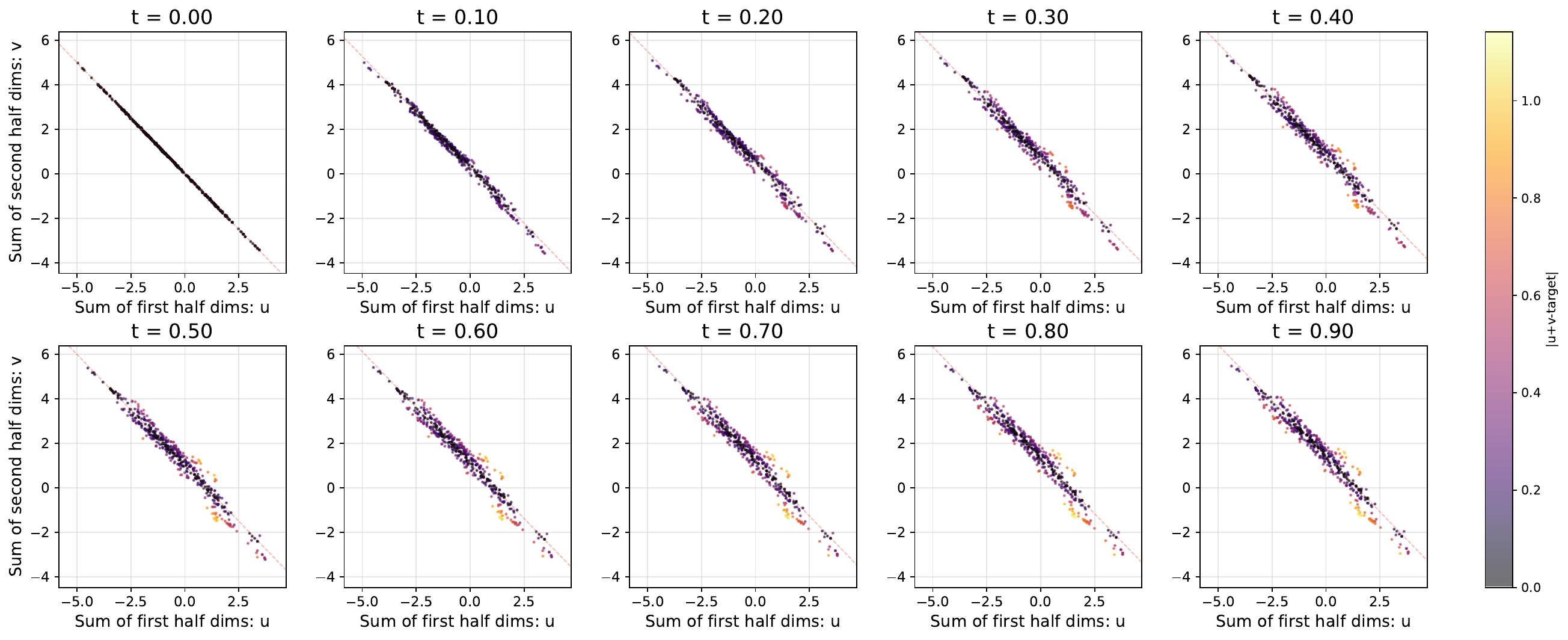}}
\caption{Projected adaptive sample trajectories generated by PMSM for the 6D variable-speed Burgers' equation \eqref{eq:dBurgers_general}. The projection uses $u=x_1+x_2+x_3$ and $v=x_4+x_5+x_6$}
\label{fig:6d_burgers_samples}
\end{figure}

\begin{figure}[htbp]
\centering
\includegraphics[width=0.88\textwidth]{\detokenize{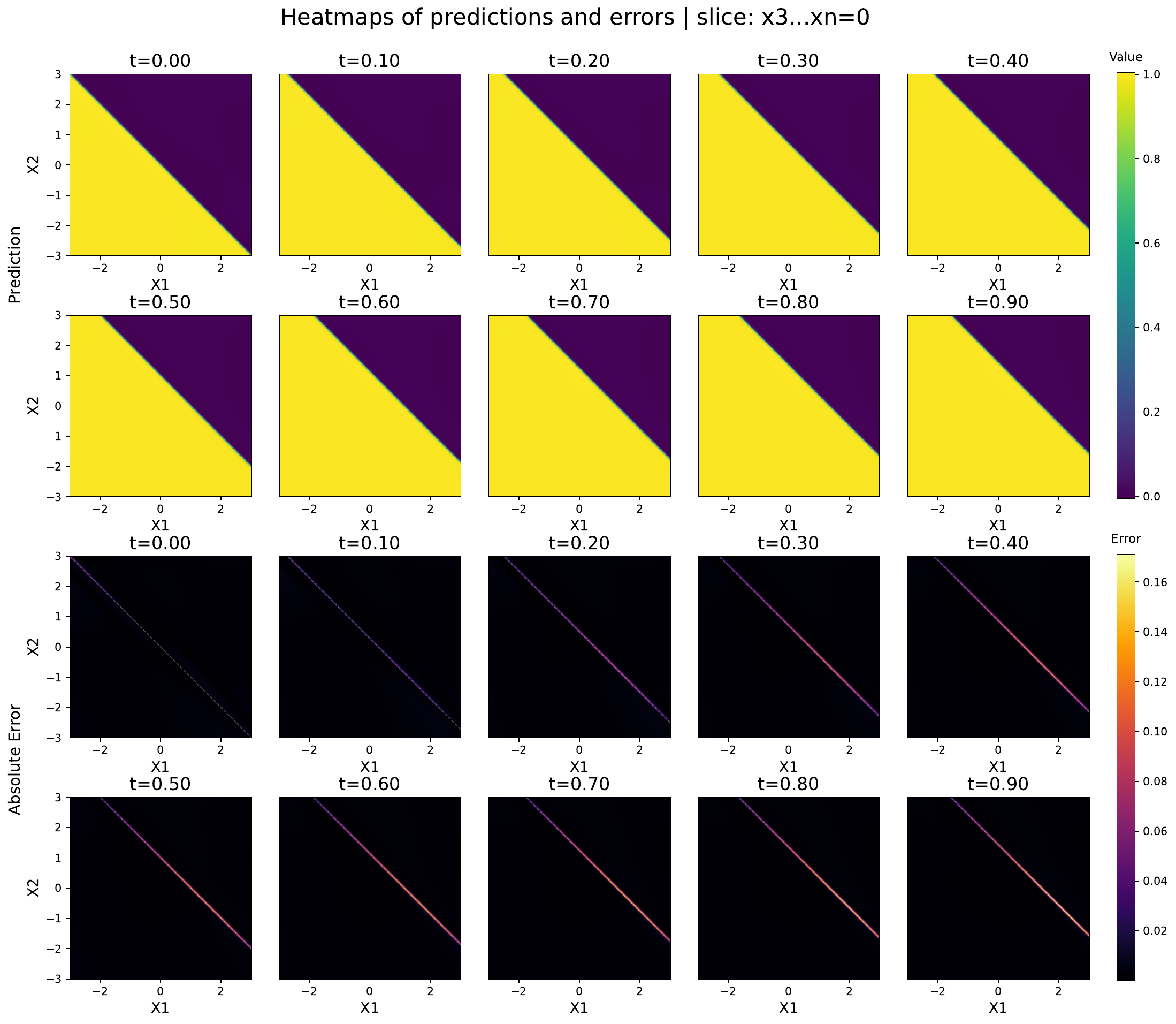}}
\caption{PINNs solution (top) and absolute error (bottom) slice heatmaps for the 6D variable-speed Burgers' equation \eqref{eq:dBurgers_general}. The visualization uses the slice $x_3=x_4=x_5=x_6=0$}
\label{fig:6d_burgers_heatmap_PINN}
\end{figure}

Since the full six-dimensional sample points cannot be visualized directly, Figure~\ref{fig:6d_burgers_heatmap} and \ref{fig:6d_burgers_heatmap_PINN} show a 2D slice of the predicted solution and the corresponding error obtained by PMSM and PINNs, where $x_3=x_4=x_5=x_6=0$. The figure again confirms that the dominant error remains concentrated near the moving front. Figure~\ref{fig:6d_burgers_samples} further shows that PMSM is able to keep the adaptive points close to the decelerating front trajectory rather than transporting them with a fixed speed. This behavior is precisely what the residual time derivative is intended to capture. MSM achieves comparable results, but PMSM exhibits better stability and solution accuracy. In this example, we demonstrate that PMSM is also capable of effectively solving high-dimensional problems.

\subsection{Effect of WR-PMSM}

We next examine the practical effect of WR-PMSM in Algorithm~\ref{alg:wr_pmsm}. The goal of this variant is to keep the extension stage training problem local in time while the adaptive trajectories continue to be rolled out over the global time axis. In all four runs reported in Table~\ref{tab:windowed_reset_results}, we use a fixed active window with $W=10$.

\begin{table}[htbp]
\caption{Comparison between PMSM and WR-PMSM. WR-PMSM uses a fixed active window with $W=10$.}
\centering
\scriptsize
\begin{tabular}{@{}p{2.35cm}p{1.55cm}p{2.15cm}p{1.95cm}p{1.55cm}@{}}
\toprule
Problem & Schedule & Relative $L^2$ & $L^\infty$ & Memory \\
\midrule
2D Burgers' & PMSM & $3.878\times 10^{-4}$ & $5.403\times 10^{-3}$ & $0.636$ GiB\\
 & WR-PMSM & $3.797\times 10^{-4}$ & $8.155\times 10^{-3}$ & $0.304$ GiB \\
\midrule
2D parabolic & PMSM & $4.338\times 10^{-2}$ & $5.618\times 10^{-2}$ & $0.643$ GiB \\
 & WR-PMSM & $8.674\times 10^{-2}$ & $8.994\times 10^{-2}$ & $0.216$ GiB \\
\midrule
3D Fokker--Planck & PMSM & $7.267\times 10^{-3}$ & $3.127\times 10^{-2}$ & $2.414$ GiB \\
 & WR-PMSM & $9.853\times 10^{-3}$ & $3.029\times 10^{-2}$ & $1.064$ GiB \\
\midrule
6D Burgers' & PMSM & $2.456\times 10^{-4}$ & $1.053\times 10^{-2}$ & $1.595$ GiB\\
 & WR-PMSM & $4.160\times 10^{-4}$ & $1.701\times 10^{-2}$ & $0.759$ GiB \\
\bottomrule
\end{tabular}
\label{tab:windowed_reset_results}
\end{table}

\begin{figure}[htbp]
\centering
\includegraphics[width=0.72\textwidth]{\detokenize{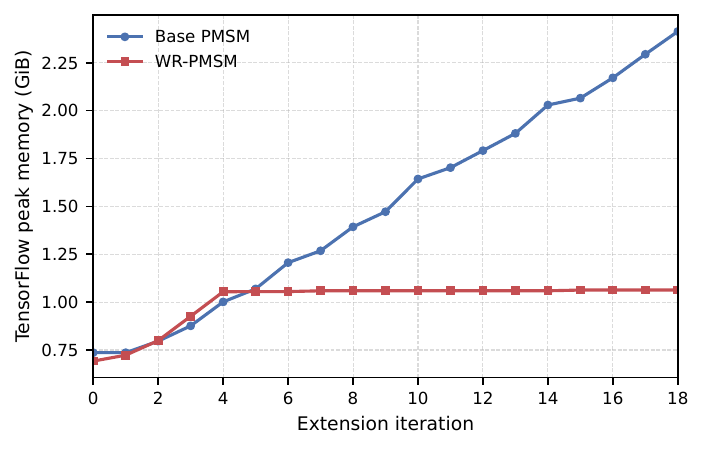}}
\caption{Extension stage TensorFlow peak memory on the 3D Fokker--Planck problem. The PMSM run uses an inactive reset setting, while WR-PMSM uses a fixed active window with $W=10$. The experimental parameter settings are as described in Subsection~\ref{subsec:3d-fokker-planck}}
\label{fig:windowed_reset_memory_trend}
\end{figure}

Table~\ref{tab:windowed_reset_results} reports the results and memory usage of PMSM and WR-PMSM across the four examples. We can observe that WR-PMSM significantly reduces computational resource consumption, while the loss in accuracy remains minor. The difference between these two methods is most relevant for long time training, where the base extension stage repeatedly revisits an increasingly long time axis. To display this effect, Figure~\ref{fig:windowed_reset_memory_trend} compares the extension stage TensorFlow peak memory on the 3D Fokker--Planck problem. Without an effective reset, the active time range grows with the rollout and the peak memory increases from $0.736$ GiB to $2.413$ GiB. With the fixed window, the same quantity stabilizes near $1.060$ GiB after the active window is filled.

Our tests confirm that this fixed-window configuration preserves solution accuracy while substantially lowering memory and computational cost for long time predictions. As a result, the total computational cost scales linearly with the prediction length, in contrast to the quadratic growth exhibited by the base algorithm without windowing.

\section{Conclusion}\label{conclusion}

This paper proposes the predictive moving sample method (PMSM) for numerically solving time-dependent partial differential equations that feature moving singular structures. PMSM is built upon the probabilistic transport framework of the moving sample method (MSM), inheriting its dual-network architecture and the core idea of residual-driven transport of the sampling distribution. To accommodate long-time prediction scenarios, two key improvements are introduced.

First, we replace the full-time-domain synchronous iteration of MSM with a stepwise progressive time-extension strategy. This strategy decouples the training scale of the extension phase from the total prediction duration, making long-time prediction computationally feasible. Second, we simplify the velocity field loss function. Under the approximation that the time interval is sufficiently small, we remove the full-domain time integration and the global normalization term from the original MSM loss, so that the velocity field training at each step relies only on the local residual information of the current time slice, thereby substantially reducing the per-step training cost. We further propose a variant, WR-PMSM, which employs an active training window strategy and periodic resets of the initial condition reference. By confining the optimization in the extension phase to a fixed-size time window, WR-PMSM allows the computational cost to grow nearly linearly with the prediction duration.

Experimental results on four representative benchmark equations demonstrate that, under a comparable budget of interior collocation points, PMSM consistently achieves lower relative $L^2$ and $L^\infty$ errors than both standard PINNs and the original MSM. WR-PMSM maintains comparable solution accuracy while substantially reducing memory consumption and computational overhead in long-time scenarios.

Several directions are worth exploring in future work. First, the method’s sensitivity to the quality of the initial sample set should be reduced. In high-dimensional problems, the initial adaptive sample set may not adequately cover difficult regions, so the effectiveness of subsequent sample transport is bounded by the information content of the initial seeds. Second, more explicit sample quality control and online rejection mechanisms can be designed, such as identifying and discarding sample points that have become invalid or have drifted away from the target region during training, preventing them from being recurrently propagated and misleading the velocity field updates. Third, a tighter coupling between adaptive sample collection and error estimation of the solution network can be investigated, for example, by adaptively determining the number of new samples added in each round according to local error estimates, instead of using a fixed per-round sampling budget.

\section*{Acknowledgment}

This work is partly supported by the Strategic Priority Research Program of the Chinese Academy of Sciences under grant number XDA0480504, by NSFC under Grant No. 12571471 and 12494543, by the Natural Science Foundation of Shandong Province (Grant No. ZR2024MA057), and by Open Foundation of the State Key Laboratory of Mathematical Sciences (Grant No. SLMS-2025-KFKT-TD-11)
J. Zhai is also partially supported by the ShanghaiTech Faculty Start-Up Fund. 

\section*{Code Availability}
The source code for this study is available at \url{https://github.com/pmsm-pinn/pmsm-pinn}

\bibliographystyle{unsrt}  
\bibliography{ref} 

@incollection{boor1973good,
  author    = {de Boor, Carl},
  title     = {Good Approximation by Splines with Variable Knots},
  booktitle = {Spline Functions and Approximation Theory},
  editor    = {Meir, A. and Sharma, A.},
  series    = {International Series of Numerical Mathematics},
  volume    = {21},
  pages     = {57--72},
  publisher = {Birkh{\"a}user},
  year      = {1973}
}

@article{white1979selection,
  title={On selection of equidistributing meshes for two-point boundary-value problems},
  author={White, Jr, Andrew B},
  journal={SIAM Journal on Numerical Analysis},
  volume={16},
  number={3},
  pages={472--502},
  year={1979},
  publisher={SIAM}
}

@article{huang1994MMPDES,
  author  = {Huang, Weizhang and Ren, Yuhe and Russell, Robert D.},
  title   = {Moving Mesh Partial Differential Equations ({MMPDEs}) Based on the Equidistribution Principle},
  journal = {SIAM Journal on Numerical Analysis},
  volume  = {31},
  number  = {3},
  pages   = {709--730},
  year    = {1994}
}

@book{huang2010adaptive,
  title={Adaptive moving mesh methods},
  author={Huang, Weizhang and Russell, Robert D},
  volume={174},
  year={2010},
  publisher={Springer Science \& Business Media}
}

@article{raissi2019physics,
  title={Physics-informed neural networks: A deep learning framework for solving forward and inverse problems involving nonlinear partial differential equations},
  author={Raissi, Maziar and Perdikaris, Paris and Karniadakis, George E},
  journal={Journal of Computational physics},
  volume={378},
  pages={686--707},
  year={2019},
  publisher={Elsevier}
}

@inproceedings{krishnapriyan2021characterizing,
  title={Characterizing possible failure modes in physics-informed neural networks},
  booktitle={Advances in Neural Information Processing Systems},
  author={Krishnapriyan, Aditi and Gholami, Amir and Zhe, Shandian and Kirby, Robert and Mahoney, Michael},
  journal={Advances in neural information processing systems},
  volume={34},
  pages={26548--26560},
  year={2021}
}

@article{wang2021gradient,
  title={Understanding and mitigating gradient flow pathologies in physics-informed neural networks},
  author={Wang, Sifan and Teng, Yujun and Perdikaris, Paris},
  journal={SIAM Journal on Scientific Computing},
  volume={43},
  number={5},
  pages={A3055--A3081},
  year={2021},
  publisher={SIAM}
}

@article{gao2023failure,
  title={Failure-informed adaptive sampling for {PINN}s},
  author={Gao, Zhiwei and Yan, Liang and Zhou, Tao},
  journal={SIAM Journal on Scientific Computing},
  volume={45},
  number={4},
  pages={A1971--A1994},
  year={2023},
  publisher={SIAM}
}

@article{gao2023active,
  title={Active learning based sampling for high-dimensional nonlinear partial differential equations},
  author={Gao, Wenhan and Wang, Chunmei},
  journal={Journal of Computational Physics},
  volume={475},
  pages={111848},
  year={2023},
  publisher={Elsevier}
}

@article{tang2023das,
  title={{DAS-PINNs}: A deep adaptive sampling method for solving high-dimensional partial differential equations},
  author={Tang, Kejun and Wan, Xiaoliang and Yang, Chao},
  journal={Journal of Computational Physics},
  volume={476},
  pages={111868},
  year={2023},
  publisher={Elsevier}
}

@inproceedings{tang2024adversarial,
  title={Adversarial adaptive sampling: Unify {PINN} and optimal transport for the approximation of PDEs},
  author={Tang, Kejun and Zhai, Jiayu and Wan, Xiaoliang and Yang, Chao},
  booktitle={International Conference on Learning Representations},
  volume={2024},
  pages={54641--54656},
  year={2024}
}

@article{bruna2024neural,
  title={Neural Galerkin schemes with active learning for high-dimensional evolution equations},
  author={Bruna, Joan and Peherstorfer, Benjamin and Vanden-Eijnden, Eric},
  journal={Journal of Computational Physics},
  volume={496},
  pages={112588},
  year={2024},
  publisher={Elsevier}
}

@article{Yuxiao_Coupling_2023,
  title={Coupling parameter and particle dynamics for adaptive sampling in Neural Galerkin schemes},
  author={Wen, Yuxiao and Vanden-Eijnden, Eric and Peherstorfer, Benjamin},
  journal={Physica D: Nonlinear Phenomena},
  volume={462},
  pages={134129},
  year={2024},
  publisher={Elsevier}
}

@book{ambrosio2008gradient,
  title={Gradient Flows: In Metric Spaces and in the Space of Probability Measures},
  author={Ambrosio, Luigi and Gigli, Nicola and Savare, Giuseppe},
  year={2008},
  publisher={Springer Science \& Business Media}
}

@article{metropolis1953equation,
  title={Equation of state calculations by fast computing machines},
  author={Metropolis, Nicholas and Rosenbluth, Arianna W and Rosenbluth, Marshall N and Teller, Augusta H and Teller, Edward},
  journal={The journal of chemical physics},
  volume={21},
  number={6},
  pages={1087--1092},
  year={1953},
  publisher={American Institute of Physics}
}

@article{hastings1970monte,
  title={{Monte Carlo} sampling methods using Markov chains and their applications},
  author={Hastings, W Keith},
  journal ={Biometrika},
  volume={57},
  number={1},
  pages={97--109},
  year={1970},
  publisher={Oxford University Press}
}

@article{duane1987hybrid,
  author  = {Duane, Simon and Kennedy, A. D. and Pendleton, Brian J. and Roweth, Duncan},
  title   = {Hybrid {Monte Carlo}},
  journal = {Physics Letters B},
  volume  = {195},
  number  = {2},
  pages   = {216--222},
  year    = {1987}
}

@article{wu2023rar,
  title={A comprehensive study of non-adaptive and residual-based adaptive sampling for physics-informed neural networks},
  author={Wu, Chenxi and Zhu, Min and Tan, Qinyang and Kartha, Yadhu and Lu, Lu},
  journal={Computer Methods in Applied Mechanics and Engineering},
  volume={403},
  pages={115671},
  year={2023},
  publisher={Elsevier}
}

@article{yang2024moving,
  title={Moving sampling physics-informed neural networks induced by moving mesh PDE},
  author={Yang, Yu and Yang, Qihong and Deng, Yangtao and He, Qiaolin},
  journal={Neural Networks},
  volume={180},
  pages={106706},
  year={2024},
  publisher={Elsevier}
}

@article{xu2026moving,
  title={Moving sample method for solving time-dependent partial differential equations},
  author={Xu, Beining and Yu, Haijun and Zhai, Jiayu and Tang, Kejun and Wan, Xiaoliang},
  journal={arXiv preprint arXiv:2601.18575},
  year={2026}
}

@incollection{arfken2011mathematical,
  title={Helmholtz’s Theorem},
  booktitle={Mathematical methods for physicists: a comprehensive guide},
  author={Arfken, George B and Weber, Hans J and Harris, Frank E},
  year={2011},
  publisher={Academic press},
  pages={95--100}
}

@article{WangHu2024NonuniformRandom,
	title = {Non-uniform random walk for adaptive sampling},
	author = {Wang, Rouhan and Hu, Dan},
	doi = {10.1016/j.jcp.2025.114160},
	journal={Journal of Computational Physics},
	year = {2025},
	volume={538},
	pages = {114160}
}

@article{LiuZDYu2026,
    title={A Gradient-Oriented Diffusion Sampling Method for Deep Partial Differential Equation Solvers},
    author={Liu, Shiqin and Zhao, Shiyu and Di, Yana and Yu, Haijun},
    year={2025},
    journal={preprint}
}

@article{dobson2021using,
  title={Using Coupling Methods to Estimate Sample Quality of Stochastic Differential Equations},
  author={Dobson, Matthew and Li, Yao and Zhai, Jiayu},
  journal={SIAM/ASA Journal on Uncertainty Quantification},
  volume={9},
  number={1},
  pages={135-162},
  year={2021}
}

\end{document}